\documentclass{article}
\usepackage[utf8]{inputenc}
\usepackage[margin=1in, includehead, includefoot]{geometry}
\usepackage{amsthm,amssymb,amsmath}
\usepackage{enumitem}
\usepackage{mathrsfs}
\usepackage{mathtools}
\usepackage{fancyhdr}
\usepackage{indentfirst}
\usepackage{graphicx}
\usepackage{subfigure}
\usepackage{placeins}
\usepackage{wrapfig}
\usepackage{url}
\usepackage[numbers]{natbib}
\usepackage[hidelinks]{hyperref}

\newtheorem{lemma}{Lemma}[section]
\newtheorem{theorem}{Theorem}[section]
\newtheorem{corollary}{Corollary}[section]
\newtheorem{observation}{Observation}[section]

\newtheorem{conjecture}{Conjecture}

\DeclarePairedDelimiter{\ceil}{\lceil}{\rceil}
\DeclarePairedDelimiter{\floor}{\lfloor}{\rfloor}

\newcommand{\toroman}[1]{\textit{\expandafter{\romannumeral #1\relax}}}

\newcommand{\cbeginproof}[0]{\par\noindent\textit{Proof.} }
\newcommand{\cendproof}[0]{ \qed\par\vspace{1em}}

\newcommand{\npcompleteproblem}[3]{ \par\vspace{0.5em}\noindent{\textbf{#1}}\newline \textbf{INSTANCE: } #2 \newline \textbf{QUESTION: } #3 \par\vspace{0.5em} }

\newcommand{\centered}[1]{\begin{tabular}{c} #1 \end{tabular}}

\title{On Error-detecting Open-locating-dominating sets}
\author{
    \small Devin C. Jean\\
    \small Computer Science Department \\
    \small Vanderbilt University\\
    \small \texttt{devin.c.jean@vanderbilt.edu}
    \and
    \small Suk J. Seo\\
    \small Computer Science Department\\
    \small Middle Tennessee State University\\
    \small \texttt{Suk.Seo@mtsu.edu}
}
\date{}

\begin{document}
\maketitle
\thispagestyle{empty}

\begin{abstract}
An open-dominating set $S$ for a graph $G$ is a subset of vertices where every vertex has a neighbor in $S$. 
An open-locating-dominating set $S$ for a graph $G$ is an open-dominating set such that each pair of distinct vertices in $G$ have distinct set of open-neighbors in $S$.
We consider a type of a fault-tolerant open-locating dominating set called error-detecting open-locating-dominating sets.
We present more results on the topic including its NP-completeness proof, extremal graphs, and a characterization of cubic graphs that permit an error-detecting open-locating-dominating set.
\end{abstract}

{\small \textbf{Keywords:} domination,  fault tolerant detection system, open-locating-dominating sets, cubic graphs, extremal graphs}
\newline\indent {\small \textbf{AMS subject classification: 05C69}}

\section{Introduction}

An open-locating-dominating set can model a type of detection system which determines the location of a possible ``intruder'' for a facility or a possible faulty processor in a network of processors \cite{old}.
A detection system is an extensively studied graphical concept which is also known as a watching system or discriminating codes \cite{watching-sys, discrim-sets}. 
Various detection systems have been defined based on the functionality of each detector in the system.
Other well-known and much studied detection systems include identifying codes \cite{karpovsky}, locating-dominating sets \cite{dom-loc-acyclic} (See Lobstein's Bibliography \cite{dombib} for a list of the articles in this field.).
In this paper, we consider a fault-tolerant variant of an open-locating dominating set called error-detecting open-locating-dominating sets.
We present more results on the topic including its NP-completeness proof, extremal graphs, and a characterization in cubic graphs.

\subsection*{Notations and definitions}

Let $G$ be a graph with vertices $V(G)$ and edges $E(G)$.
The \emph{open-neighborhood} of a vertex $v \in V(G)$, denoted $N(v)$, is the set of all vertices adjacent to $v$: $\{w \in V(G) : vw \in E(G)\}$.
The \emph{closed-neighborhood} of a vertex $v \in V(G)$, denoted $N[v]$, is the set of all vertices adjacent to $v$, as well as $v$ itself: $N(v) \cup \{v\}$.
An \textit{open-dominating set} (also called total-dominating) $S$ of a graph $G$ is a subset of vertices where every vertex has a neighbor in $S$. 
An \textit{open-locating-dominating set} (OLD set) $S$ of a graph $G$ is an open-dominating set such that each pair of distinct vertices in $G$ have distinct open-neighborhoods within $S$.
For an OLD set $S \subseteq V(G)$ and $u \in V(G)$, we let $N_S(u) = N(u) \cap S$ denote the \emph{dominators} of $u$ and $dom(u) = |N_S(u)|$ denote the (open) \emph{domination number} of $u$.
A vertex $v \in V(G)$ is \emph{$k$-open-dominated} by an open-dominating set $S$ if $|N_S(v)| = k$.
If $S$ is a open-dominating set and $u,v \in V(G)$, $u$ and $v$ are \emph{$k$-distinguished} if $|N_S(u) \triangle N_S(v)| \ge k$, where $\triangle$ denotes the symmetric difference.
If $S$ is a open-dominating set and $u,v \in V(G)$, $u$ and $v$ are \emph{$k^\#$-distinguished} if $|N_S(u) - N_S(v)| \ge k$ or $|N_S(v) - N_S(u)| \ge k$.
We will also use terms such as ``at least $k$-dominated'' to denote $j$-dominated for some $j \ge k$.

There are several fault-tolerant variants of OLD sets.
For example, a redundant open-locating-dominating set is resilient to a detector being destroyed or going offline \cite{ftold}.
Thus, an open-dominating set $S \subseteq V(G)$ is called a \emph{redundant open-locating-dominating (RED:OLD)} set if $\forall v \in S$, $S - \{v\}$ is an OLD set.
The focus of this paper is another variant of an OLD set  called an \textit{error-detecting open-locating-dominating} (DET:OLD) set, which is capable of correctly identifying an intruder even when at most one sensor or detector incorrectly reports that there is no intruder. 
Hence, DET:OLD sets allow for uniquely locating an intruder in a way which is resilient to up to one false negative.
The following Theorem characterizes OLD, RED:OLD, and DET:OLD sets and they are useful in constructing those sets or verifying whether a given set meets their requirements.

\begin{theorem}\label{theo:multi-param} 
An open-dominating set is
\begin{enumerate}[noitemsep, label=\roman*.]
    \item \cite{ourtri} an OLD set if and only if every pair of vertices is 1-distinguished.
    \item \cite{ftold} a RED:OLD set if and only if all vertices are at least 2-dominated and all pairs are 2-distinguished.
    \item \cite{ftold} a DET:OLD set if and only if all vertices are at least 2-dominated and all pairs are $2^\#{}$-distinguished.
\end{enumerate}
\label{theo:old-redold-detold}
\end{theorem}

\begin{figure}[ht]
    \centering
    \begin{tabular}{c@{\hspace{4em}}c@{\hspace{4em}}c}
        \includegraphics[width=0.175\textwidth]{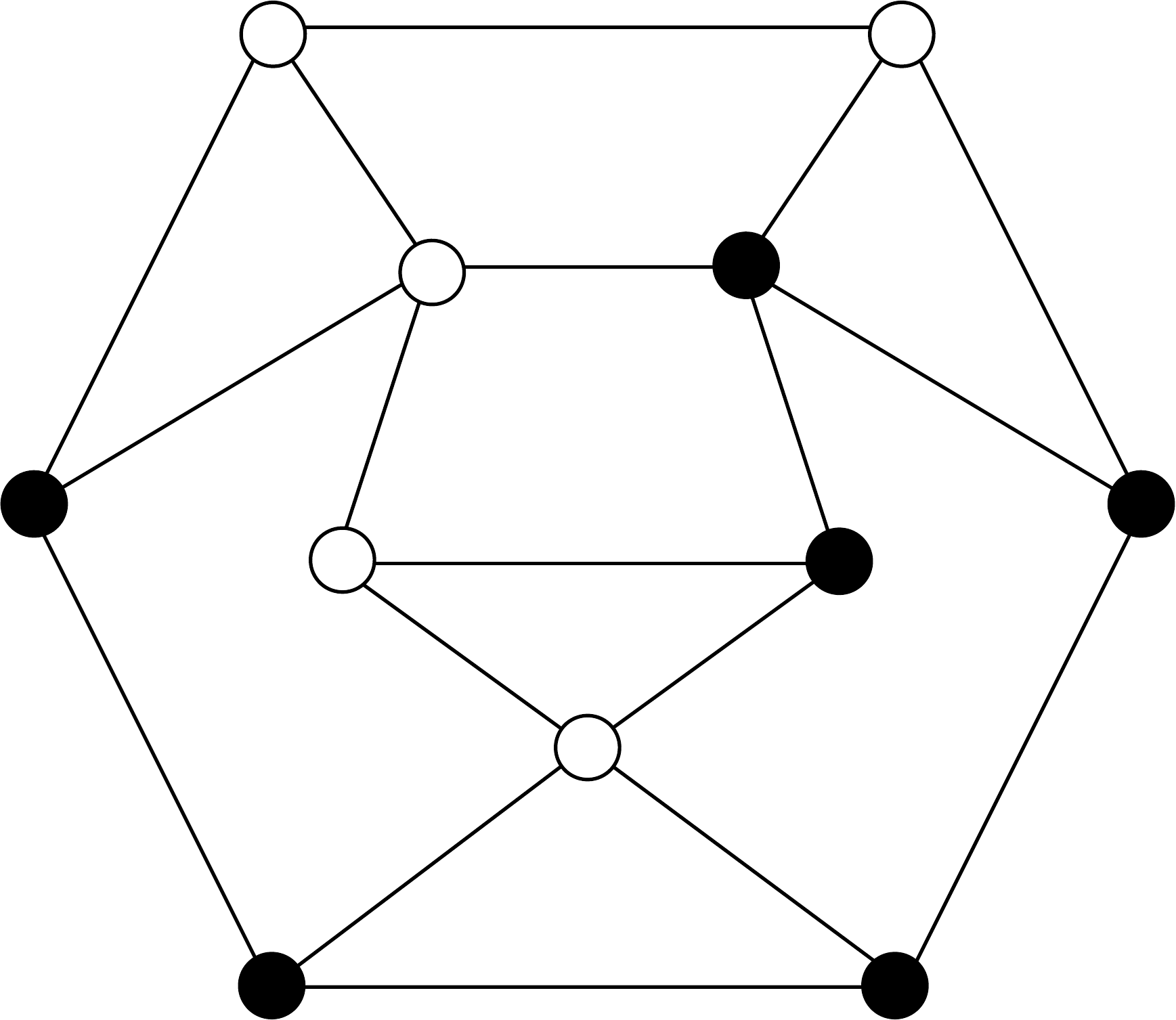} & \includegraphics[width=0.175\textwidth]{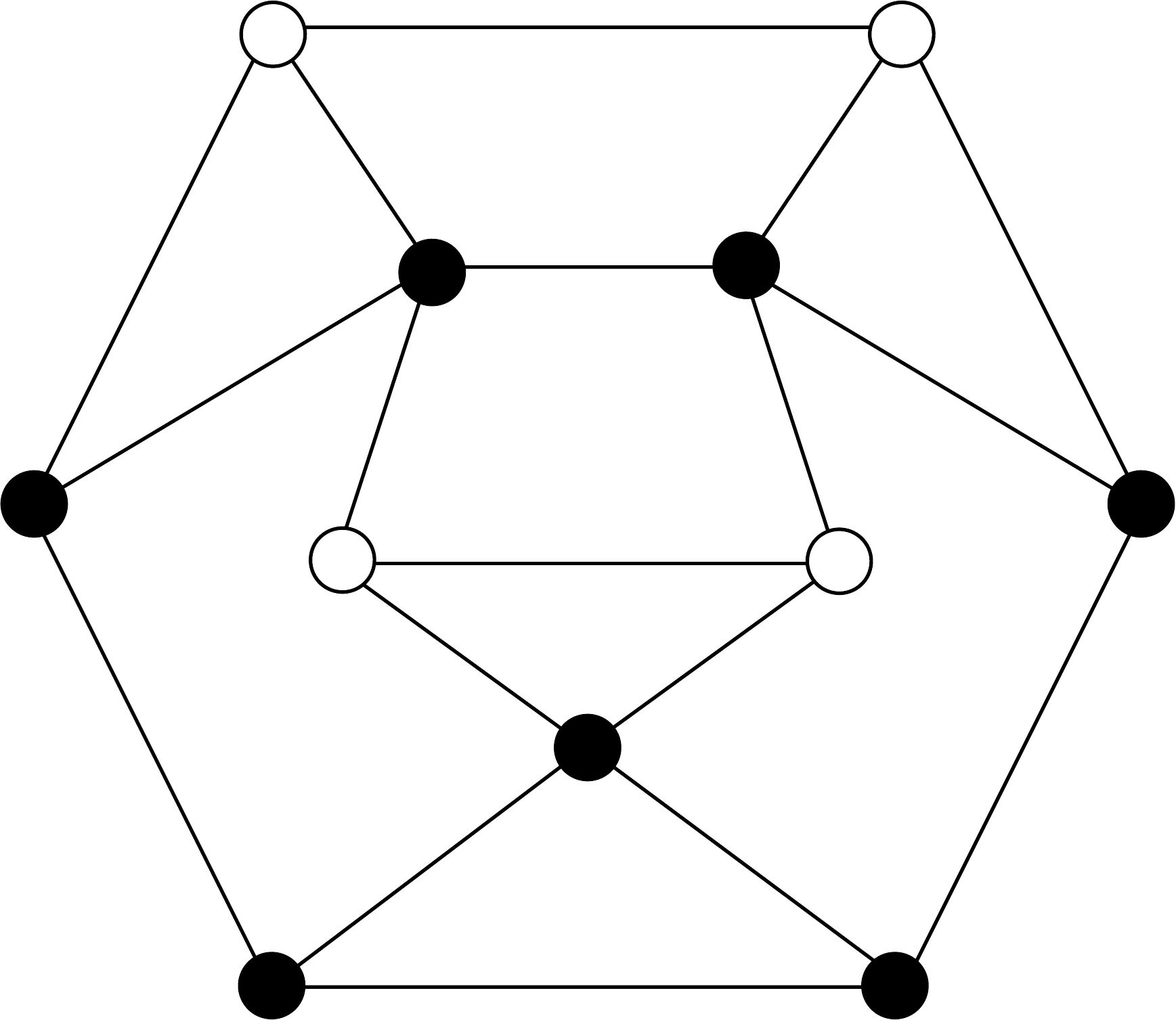} & \includegraphics[width=0.175\textwidth]{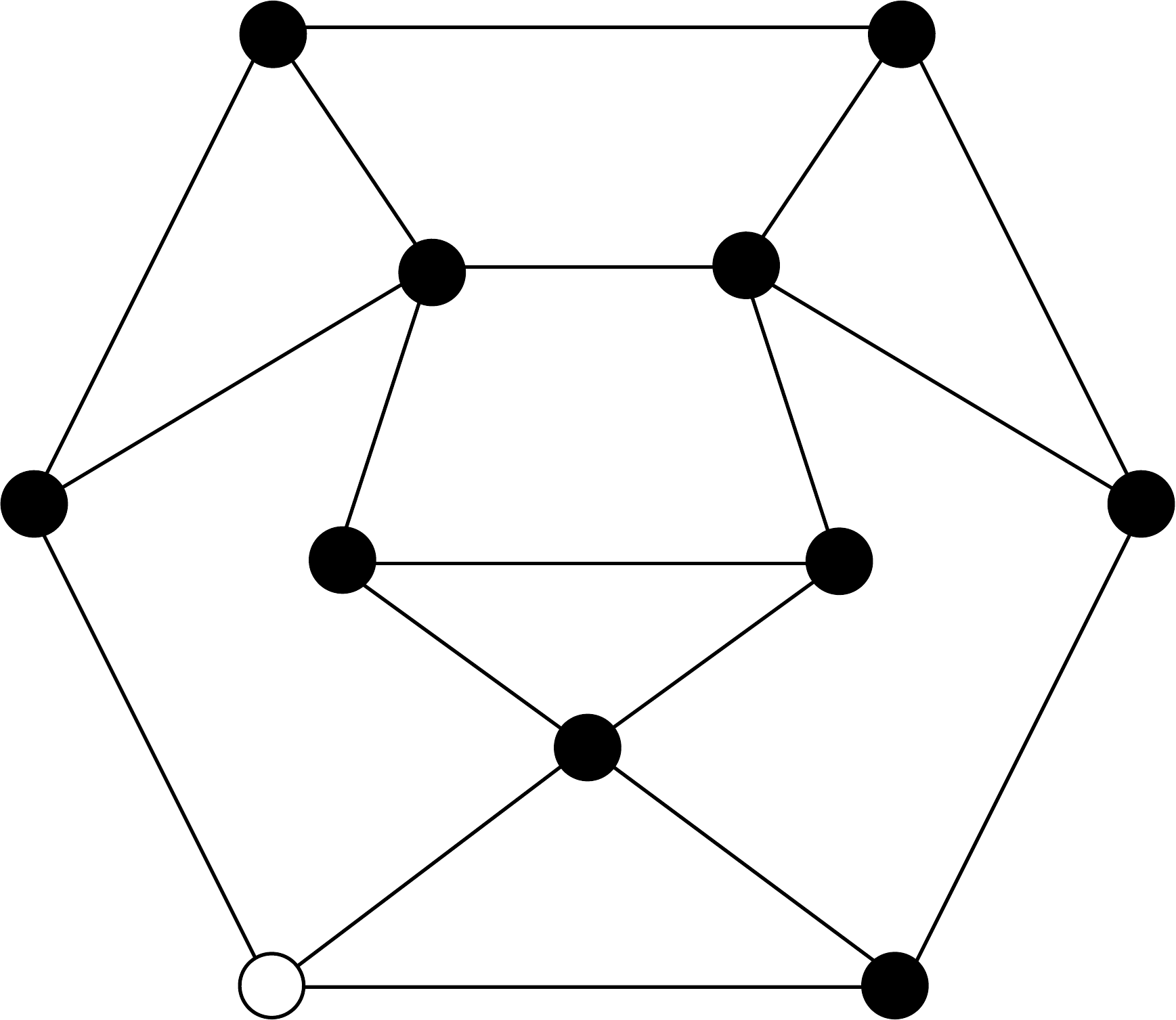} \\
        (a) & (b) & (c)
    \end{tabular}
    \caption{Optimal OLD (a), RED:OLD (b), and DET:OLD (c) sets. Shaded vertices represent detectors.}
    \label{fig:ex-finite}
\end{figure}

Naturally, our goal is to install a minimum number of detectors in any detection system. 
For finite graphs, the notations $\textrm{OLD}(G)$, $\textrm{RED:OLD}(G)$, and $\textrm{DET:OLD}(G)$ represent the cardinality of the smallest possible OLD, RED:OLD, and DET:OLD sets on graph $G$, respectively  \cite{ourtri, ftold, ft-old-cubic}.
Figure~\ref{fig:ex-finite} shows an OLD, RED:OLD, and DET:OLD sets on the given graph $G$.
We can verify that those sets of detectors meet the requirements specified in Theorem ~\ref{theo:old-redold-detold}.
There are no other sets that use fewer number of detectors for each of the three parameters, so we get $\textrm{OLD}(G) = 6$, $\textrm{RED:OLD}(G) = 7$, and $\textrm{DET:OLD}(G) = 10$.
For infinite graphs, instead of the cardinality, we measure via the \emph{density} of the subset, which is defined as the ratio of the number of detectors to the total number of vertices.
The notations $\textrm{OLD}\%(G)$, $\textrm{RED:OLD}\%(G)$, and $\textrm{DET:OLD}\%(G)$ represent the minimum density of such a set on $G$.
Note that the notion of density is also defined for finite graphs.

In Section~\ref{sec:npc}, we present a proof that the problem of determining $\textrm{DET:OLD}(G)$ for an arbitrary graph is NP-complete.
Section~\ref{sec:detold-char} shows extremal graphs with the highest density.
Sections~\ref{sec:grids} and \ref{sec:cubic} discusses DET:OLD sets in the infinite regular graphs and cubic graphs, respectively.

\section{Extremal Graphs on DET:OLD sets}\label{sec:detold-char}

In this section, we consider extremal graphs with $\textrm{DET:OLD}(G) = n$.
Let $S \subseteq V(G)$ be a DET:OLD set for $G$; because non-detectors do not aid in locating intruders, it must be that $S$ is a DET:OLD set for $G - (V(G) - S)$, implying the smallest graph with DET:OLD will have $\textrm{DET:OLD}(G) = n$.
The following theorem shows that the smallest graphs with DET:OLD have $n = 7$.

\begin{theorem}
If $G$ has a DET:OLD, then $n \ge 7$.
\end{theorem}
\begin{proof}
Assume to the contrary,  $n \le 6$.
Clearly $G$ has a cycle because DET:OLD requires $\delta(G) \ge 2$.
Firstly, we consider the cases when the smallest cycle in the graph is $C_n$, with n = 6,5,4.
If the smallest cycle is a $C_6$ subgraph $abcdef$, then $a$ and $c$ cannot be distinguished without having $n \ge 7$, a contradiction.
Suppose the smallest cycle is a $C_5$ subgraph $abcde$.
To distinguish $a$ and $c$, by symmetry we can assume that $\exists u \in N(a) - N(c) - \{a,b,c,d,e\}$.
Similarly, to distinguish $b$ and $e$ we can assume by symmetry that $\exists v \in N(b) - N(e) - \{a,b,c,d,e\}$.
If $u = v$, there would be a smaller cycle, implying $n \ge 7$, a contradiction; thus, we can assume $u \neq v$.
Suppose the smallest cycle is a $C_4$ subgraph $abcd$.
To distinguish $a$ and $c$, we can assume $p,q \in N(a) - N(c) - \{a,b,c,d\}$.
Similarly, to distinguish $b$ and $d$ we can assume by symmetry that $\exists x,y \in N(b) - N(d) - \{a,b,c,d\}$.
We know that $\{p,q\} \cap \{x,y\} = \varnothing$ because otherwise we create a smaller cycle; thus, $n \ge 8$, a contradiction.
Otherwise, we can assume that $G$ has a triangle.

Next, we will show that if $G$ contains a $K_4$ subgraph, then $n \ge 7$; let $abcd$ be the vertices of a $K_4$ subgraph in $G$.
To distinguish pairs of vertices in $abcd$, without loss of generality we can assume that $\exists x \in N(a)$, $\exists y \in N(b)$, and $\exists z \in N(c)$.
Clearly $\{x,y,z\} \cap \{a,b,c,d\} = \varnothing$ because $x$, $y$, and $z$ are used to distinguish the vertices $a$, $b$, $c$, and $d$; however, we do not yet know if $x$, $y$, and $z$ are distinct.
If $x=y=z$, then distinguishing $a$, $b$, and $c$ requires at least another two vertices, so we would have at least 7 vertices and would be done; otherwise, without loss of generality we can assume $x \neq z$.
Suppose $x=y$; we see that $a$ and $b$ are not distinguished and $n=6$, so without loss of generality let $bz \in E(G)$ to distinguish $a$ and $b$.
If $\{dx,dz\} \cap E(G) = \varnothing$, then $(d,x)$ and $(d,z)$ cannot be distinguished; otherwise without loss of generality assume $dx \in E(G)$.
We now see that $a$ and $d$ are not distinguished, but they are symmetric, so without loss of generality let $dz \in E(G)$ to distinguish $a$ and $d$.
We see that $d$ and $b$ become closed twins and cannot be distinguished.
Otherwise, we can assume $x \neq y$ and by symmetry $y \neq z$.
Thus, $n \ge 7$, and we would be done.
Thus, if $G$ has a $K_4$ subgraph, then we would be done.

Next, we will show that the existence of a ``diamond'' subgraph, which is an (almost) $K_4$ subgraph minus one edge, implies that $n \ge 7$.
Let $abcd$ be a $C_4$ subgraph and assume $ac \in E(G)$ but $cd \notin E(G)$, which forms said diamond subgraph.
To distinguish $b$ and $d$, we can assume by symmetry that $\exists u,v \in N(b) - N(d) - \{a,b,c,d\}$.
Similarly, to distinguish $a$ and $c$, we can assume that $\exists w \in N(a) - N(c) - \{a,b,c,d\}$.
We know that $u \neq v$ by assumption, so $n \le 6$ requires by symmetry that $w = u$.
We can assume that $uc \notin E(G)$ because otherwise this would create a $K_4$ subgraph and fall into a previous case.
We see that $cv \in E(G)$ is required to distinguish $u$ and $c$.
To distinguish $u$ and $v$, by symmetry we can assume $uv \in E(G)$.
Now, we see that $a$ and $v$ cannot be distinguished without creating a $K_4$ subgraph, so we are done with the diamond case.

For the final case, we know there must be a triangle, $abc$, but there cannot be any $K_4$ or diamond subgraphs.
To distinguish vertices in $abc$, we can assume by symmetry that $\exists u \in N(a) - \{a,b,c\}$ and $\exists v \in N(b) - \{a,b,c\}$; further, we know that $u \neq v$ because otherwise this would create a diamond subgraph.
We know that $\{av, bu, cu, cv\} \cap E(G) = \varnothing$ because any of these edges would create a diamond subgraph.
Suppose $uv \in E(G)$; then distinguishing $a$ and $v$ within the bounds of $n \le 6$ requires $\exists w \in N(a) - N(v) - \{a,b,c,u,v\}$.
Similarly, distinguishing $u$ and $b$ requires $w \in N(b)$; however, this creates a diamond subgraph, so we would be done.
Otherwise, we can assume that $uv \notin E(G)$.
To 2-dominate $u$ and $v$, we require $\exists p \in N(u) - \{a,b,c,u,v\}$ and $\exists q \in N(v) - \{a,b,c,u,v\}$.
However, $n \le 6$ requires that $p = q$.
We see that $u$ and $v$ cannot be distinguished without creating a diamond subgraph or having $n \ge 7$, completing the proof.
\end{proof}

\begin{figure}[ht]
    \centering
    \includegraphics[width=0.1\textwidth]{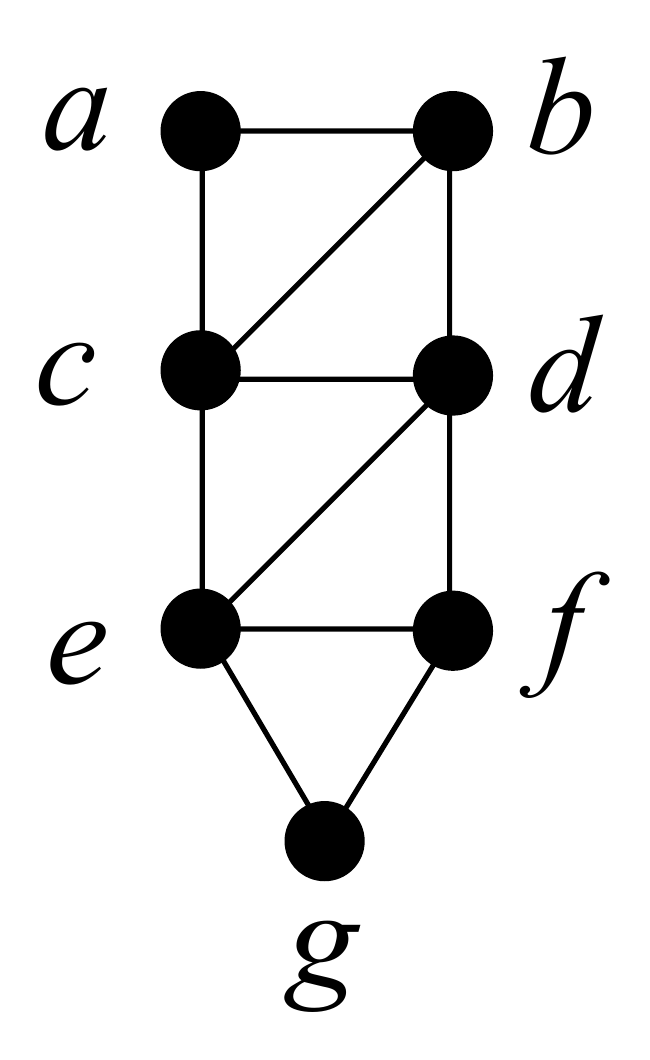}
    \caption{$G_{7,11}$ with $\textrm{DET:OLD}(G_{7,11}) = 7$}
    \label{fig:g7}
\end{figure}

Let $G_{n,m}$ have $n$ vertices and $m$ edges.
Then $G_{7,11}$, as shown in Figure~\ref{fig:g7}, is the first graph that permits a DET:OLD set in the lexicographic ordering of $(n,m)$ tuples; i.e., the graph with the smallest number of edges given the smallest number of vertices.

Next we consider extremal graphs with $\textrm{DET:OLD}(G) = n$ with the fewest number of edges.

\begin{observation}\label{obs:detold-deg2-neighbor}
If $S$ is a DET:OLD set on $G$ and $v \in V(G)$, then $v$ has at most one degree 2 neighbor.
\end{observation}

\begin{theorem}\label{theo:min-edges}
If $G_{n,m}$ has DET:OLD then $m \ge \ceil*{\frac{3n - \floor{\frac{n}{2}}}{2}}$.
\end{theorem}
\begin{proof}
Because $G$ has DET:OLD, we know $\delta(G) \ge 2$; let $p$ be the number of degree 2 vertices in $G$.
By Observation~\ref{obs:detold-deg2-neighbor}, every degree 2 vertex, $v$, must have at least one neighbor, $u$, of at least degree 3, and said neighbor $u$ is not adjacent to any degree 2 vertices other than $v$.
Thus, we can pair each of the $p$ degree 2 vertices with a unique degree 3 or higher vertex.
From this we know that $p \le \floor{\frac{n}{2}}$, and the $n-2p$ vertices that are not pairs are all at least degree 3.
Thus, $\sum_{v \in V(G)}{deg(v)} \ge (2 + 3)p + 3(n - 2p) = 3n - p \ge 3n - \floor{\frac{n}{2}}$.
However, we also know that the degree sum of any graph must be even, so we can strengthen this to $\sum_{v \in V(G)}{deg(v)} \ge 2 \ceil*{\frac{3n - \floor{\frac{n}{2}}}{2}}$.
Dividing the degree sum by 2 completes the proof.
\end{proof}

\begin{figure}[ht]
    \centering
    \begin{tabular}{c@{\hspace{4em}}c@{\hspace{4em}}c@{\hspace{4em}}c}
        \includegraphics[height=0.125\textwidth]{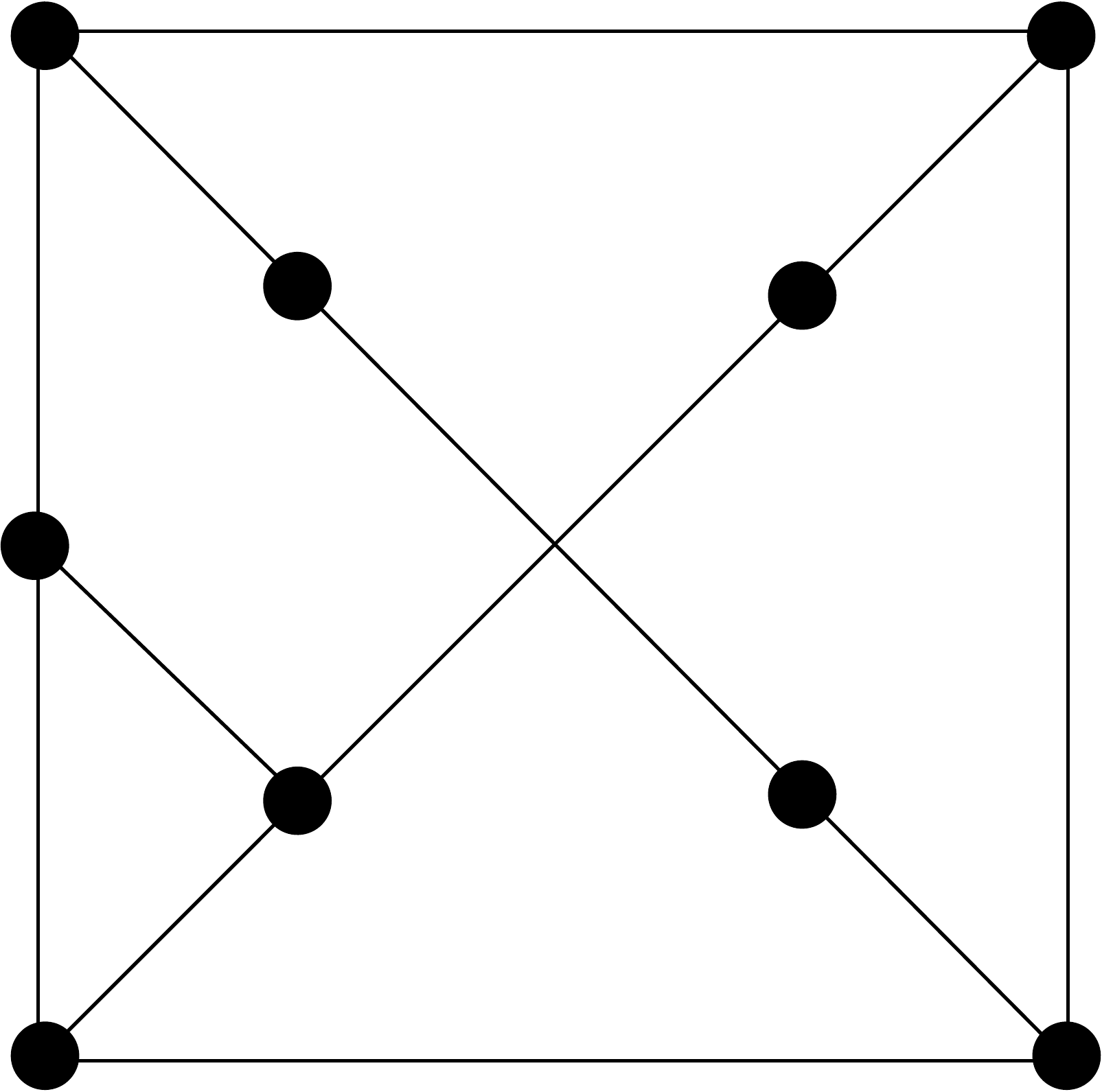} & \includegraphics[height=0.125\textwidth]{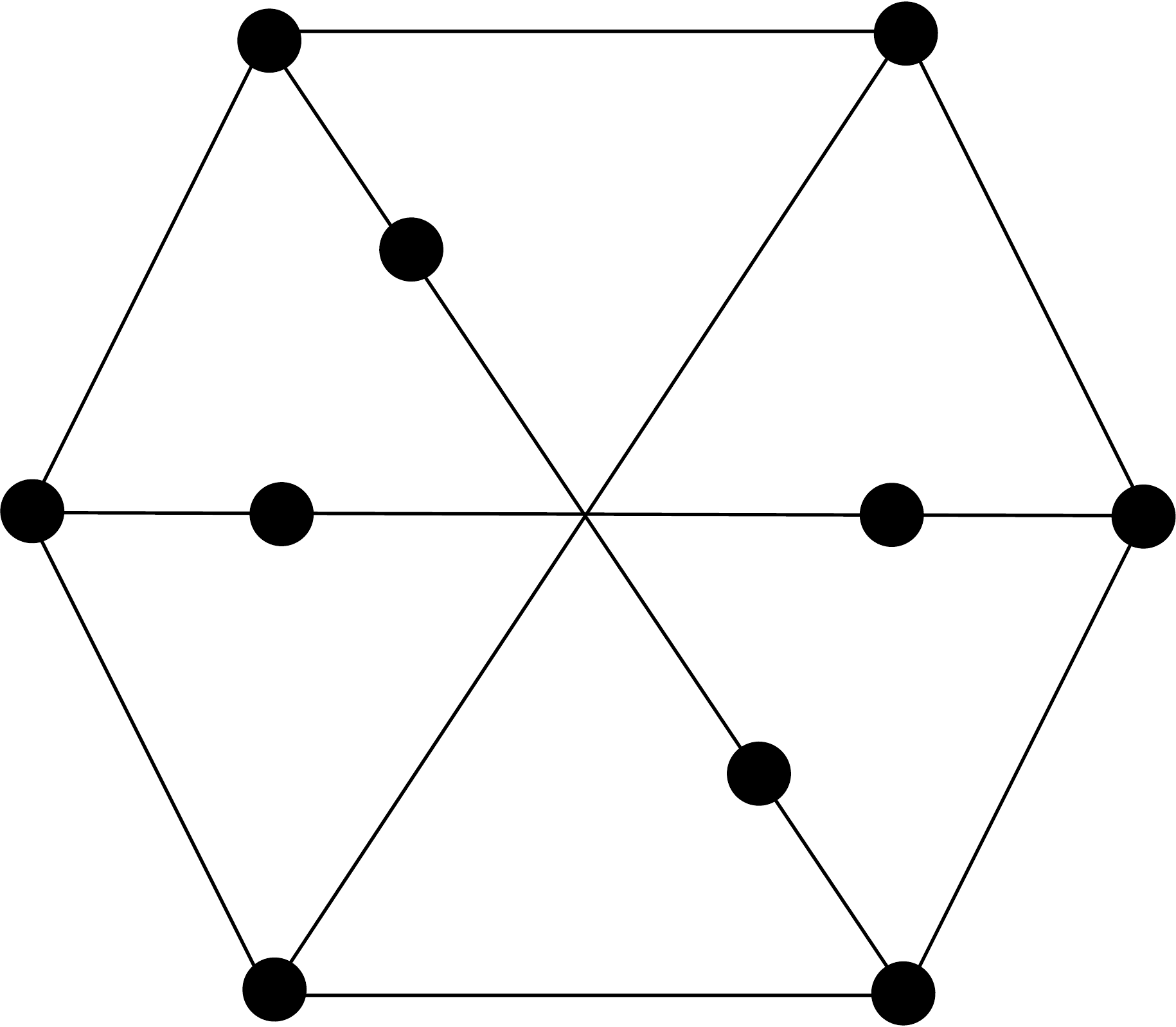} & \includegraphics[height=0.125\textwidth]{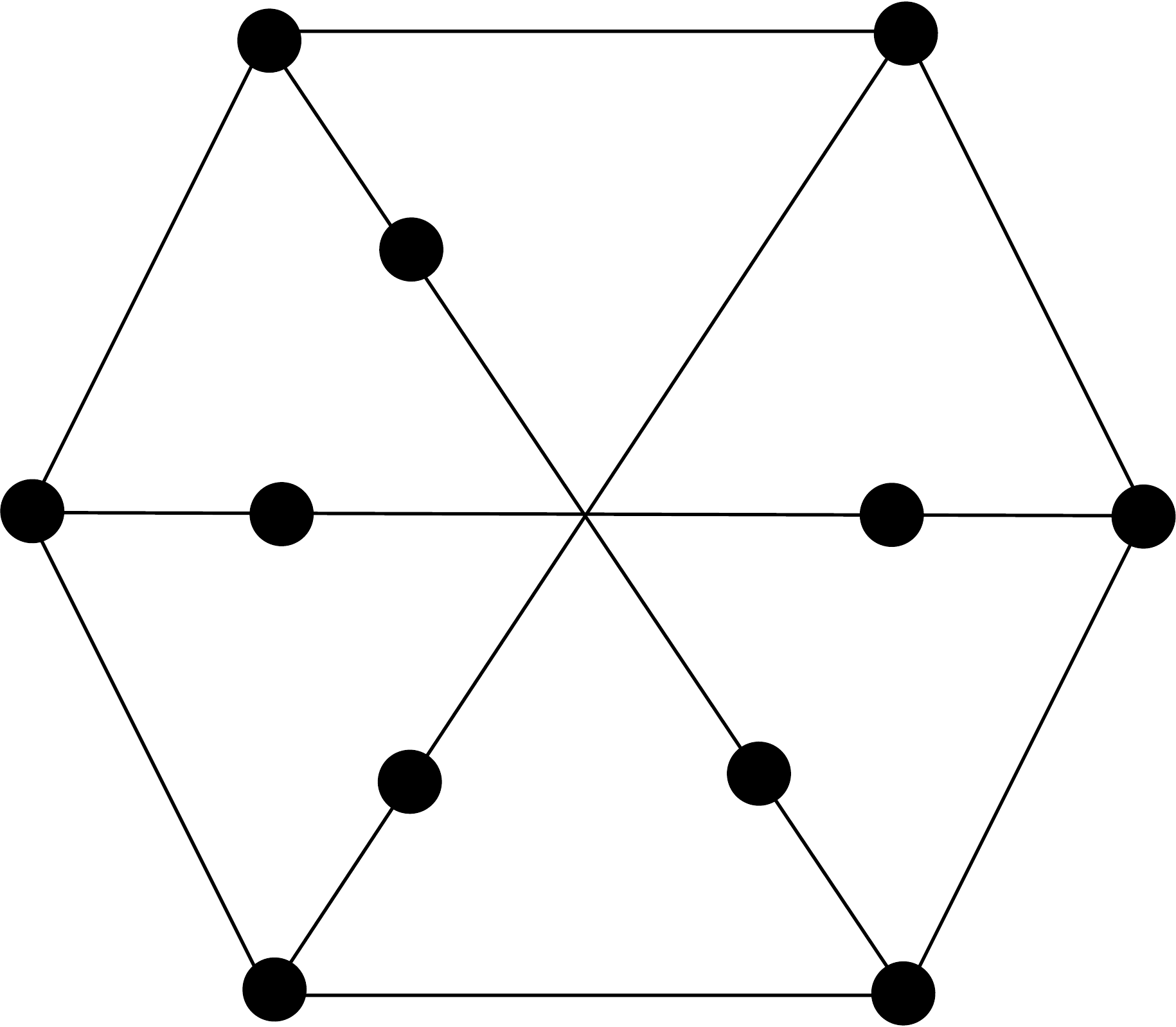} & \includegraphics[height=0.125\textwidth]{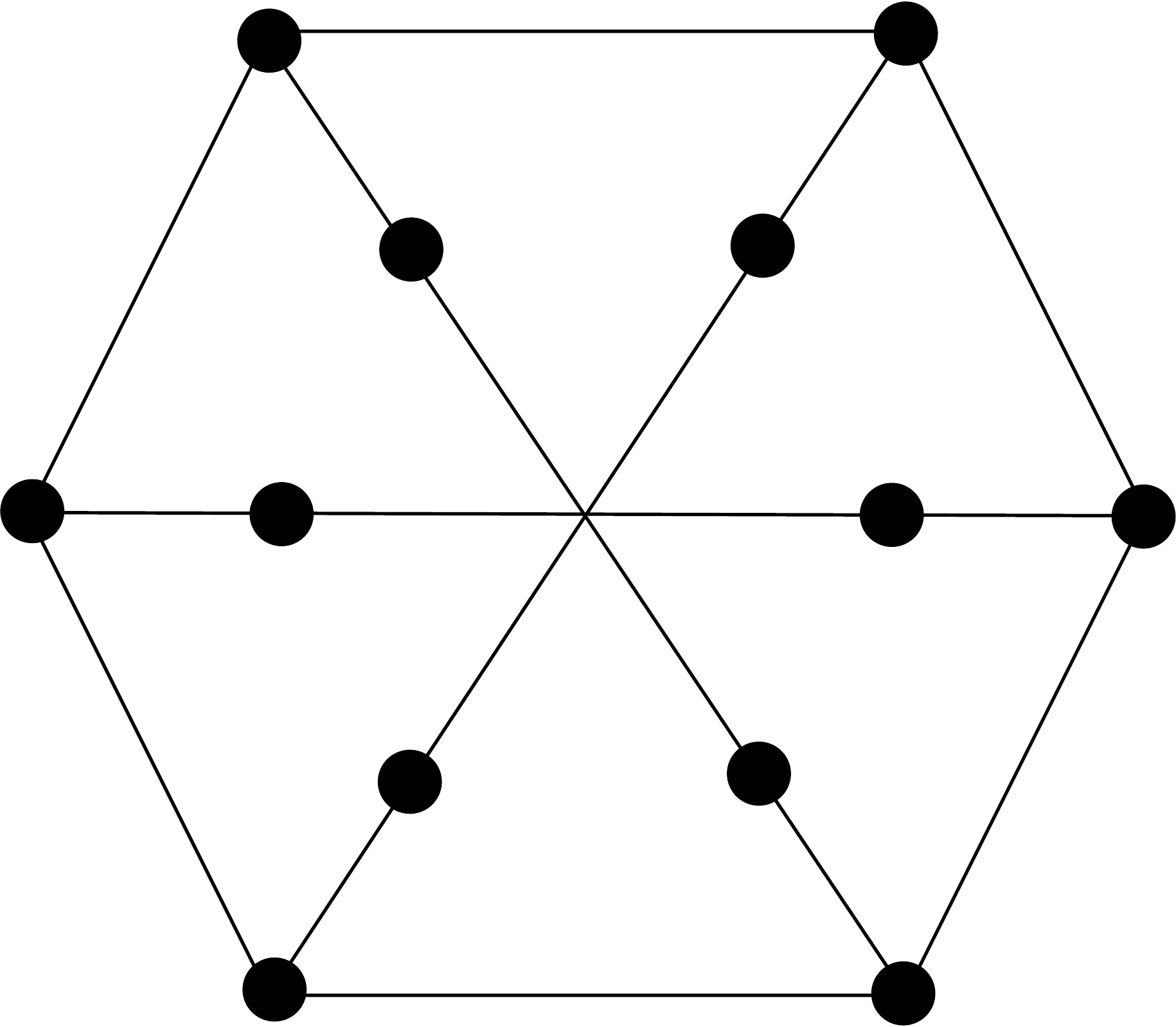} \\
        $n=9$ & $n=10$ & $n=11$ & $n=12$ \\ \\
        \includegraphics[height=0.125\textwidth]{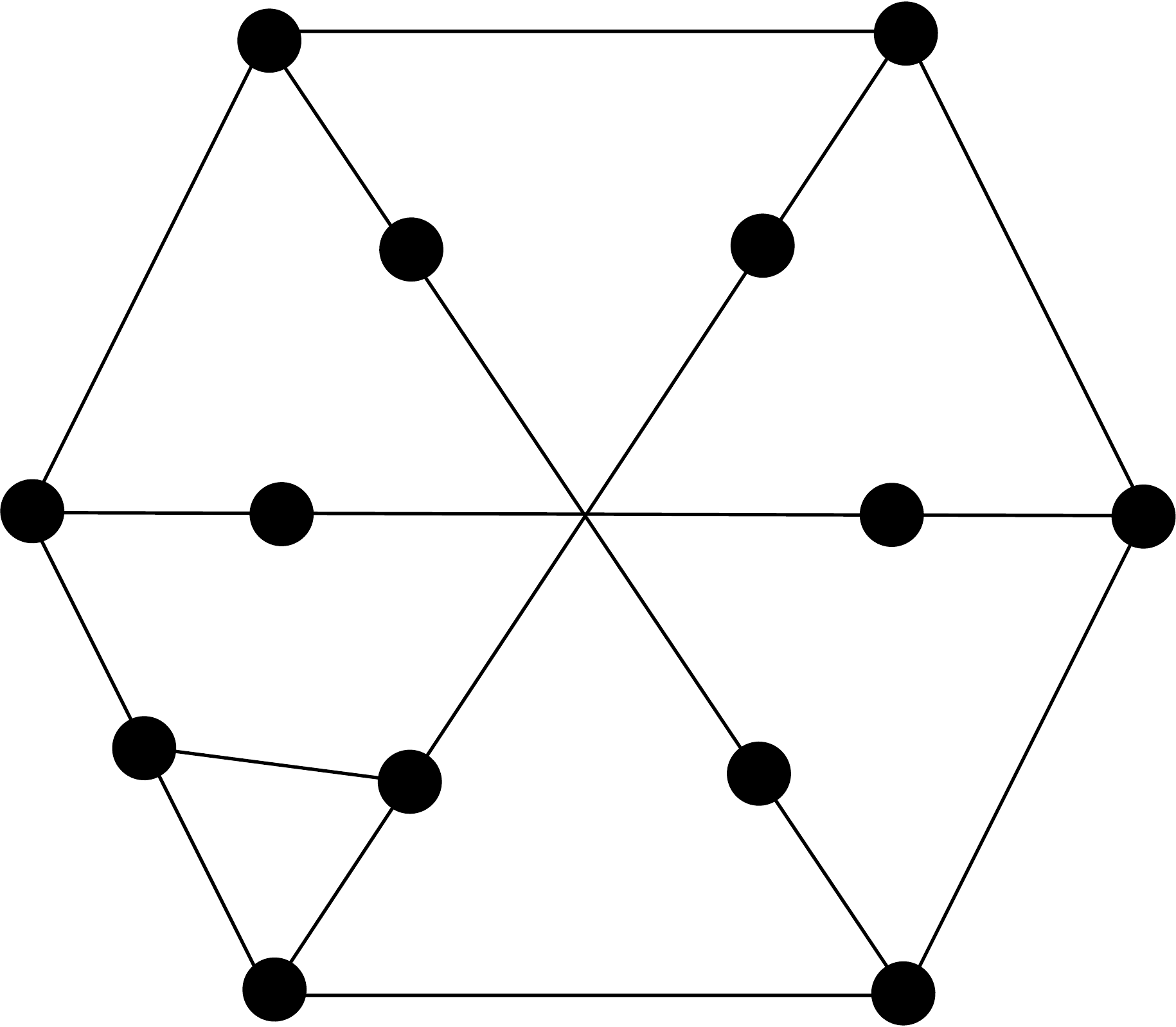} & \includegraphics[height=0.125\textwidth]{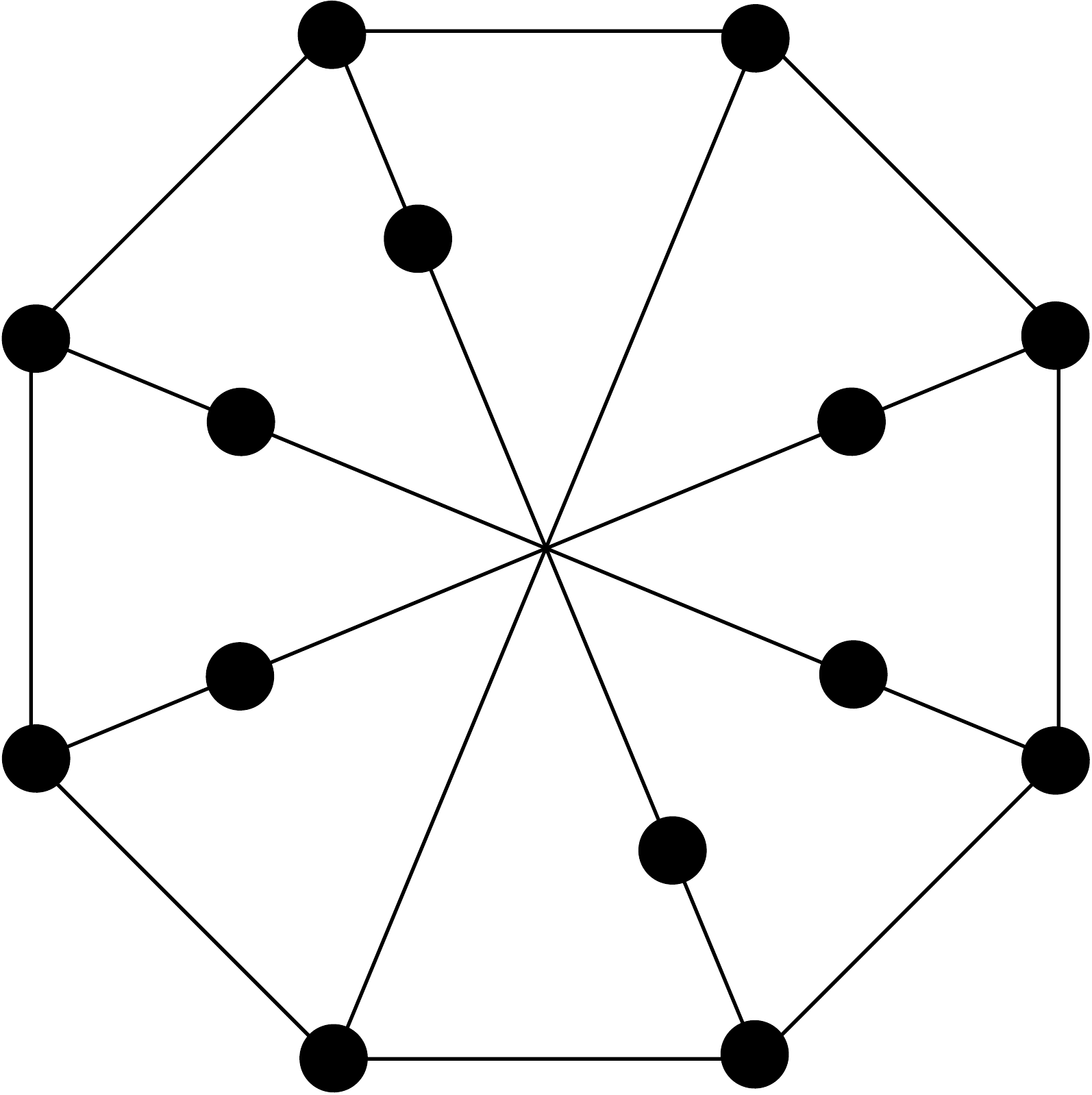} & \includegraphics[height=0.125\textwidth]{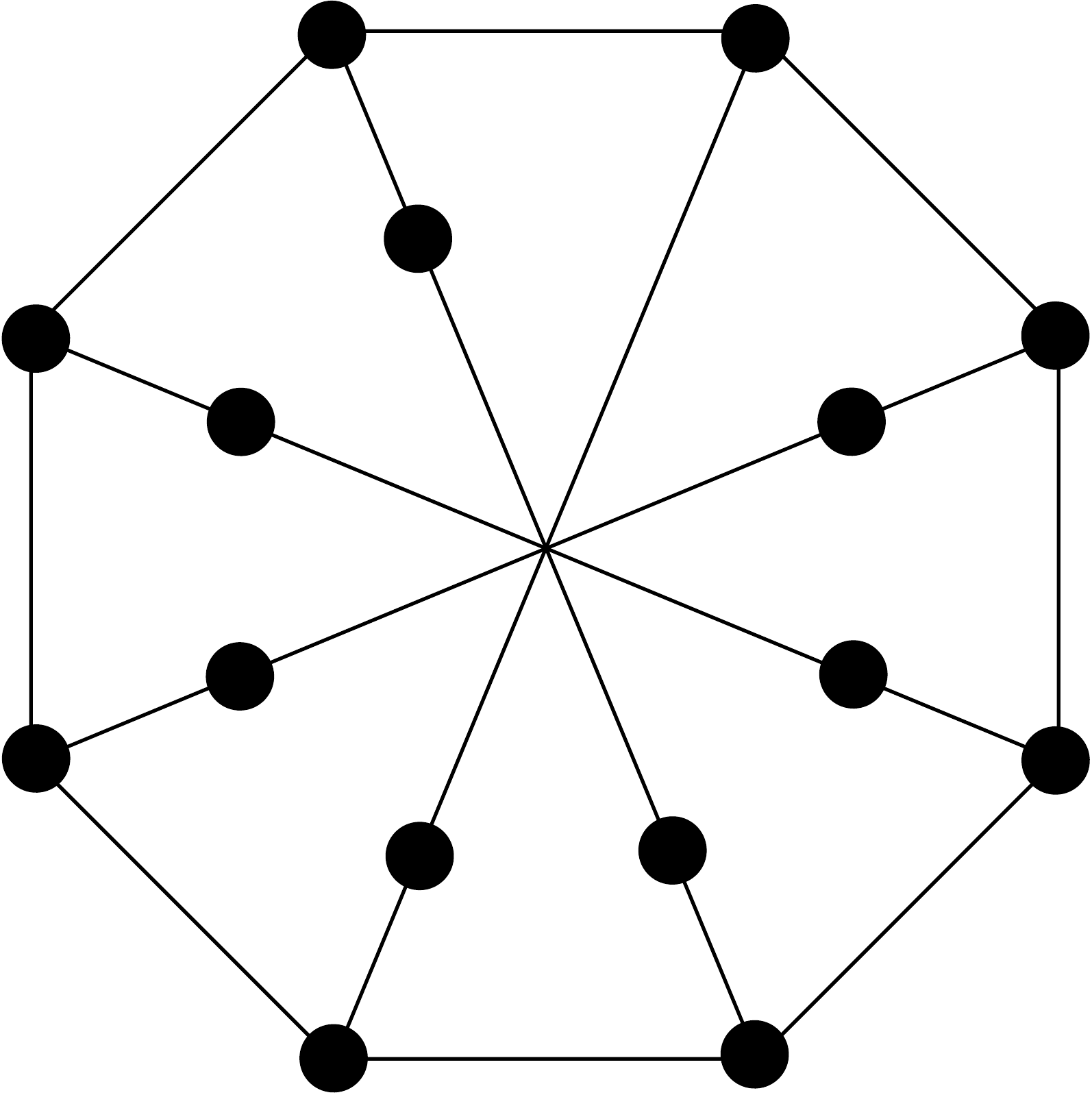} & \includegraphics[height=0.125\textwidth]{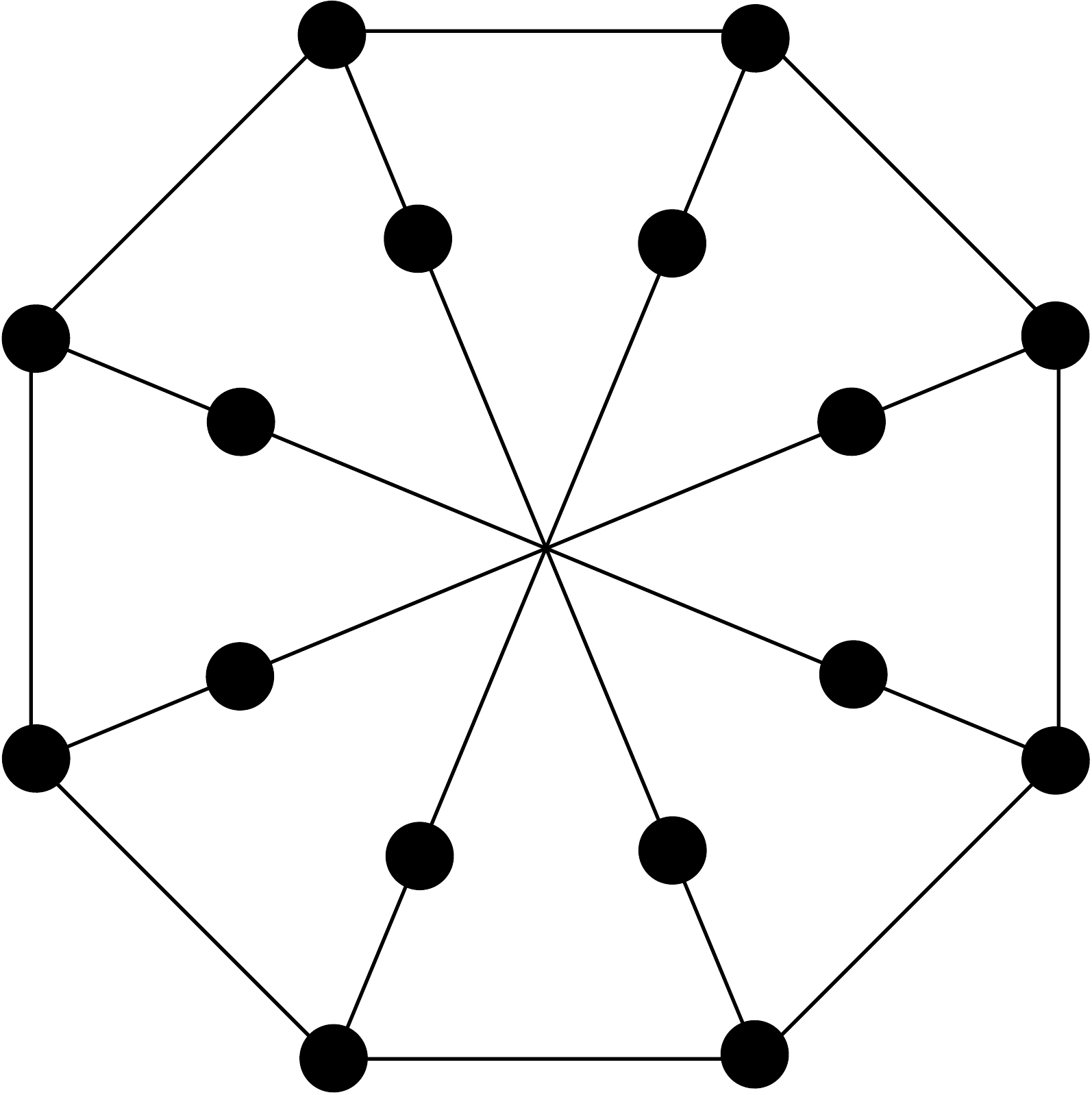} \\
        $n=13$ & $n=14$ & $n=15$ & $n=16$ \\ \\
        \includegraphics[height=0.125\textwidth]{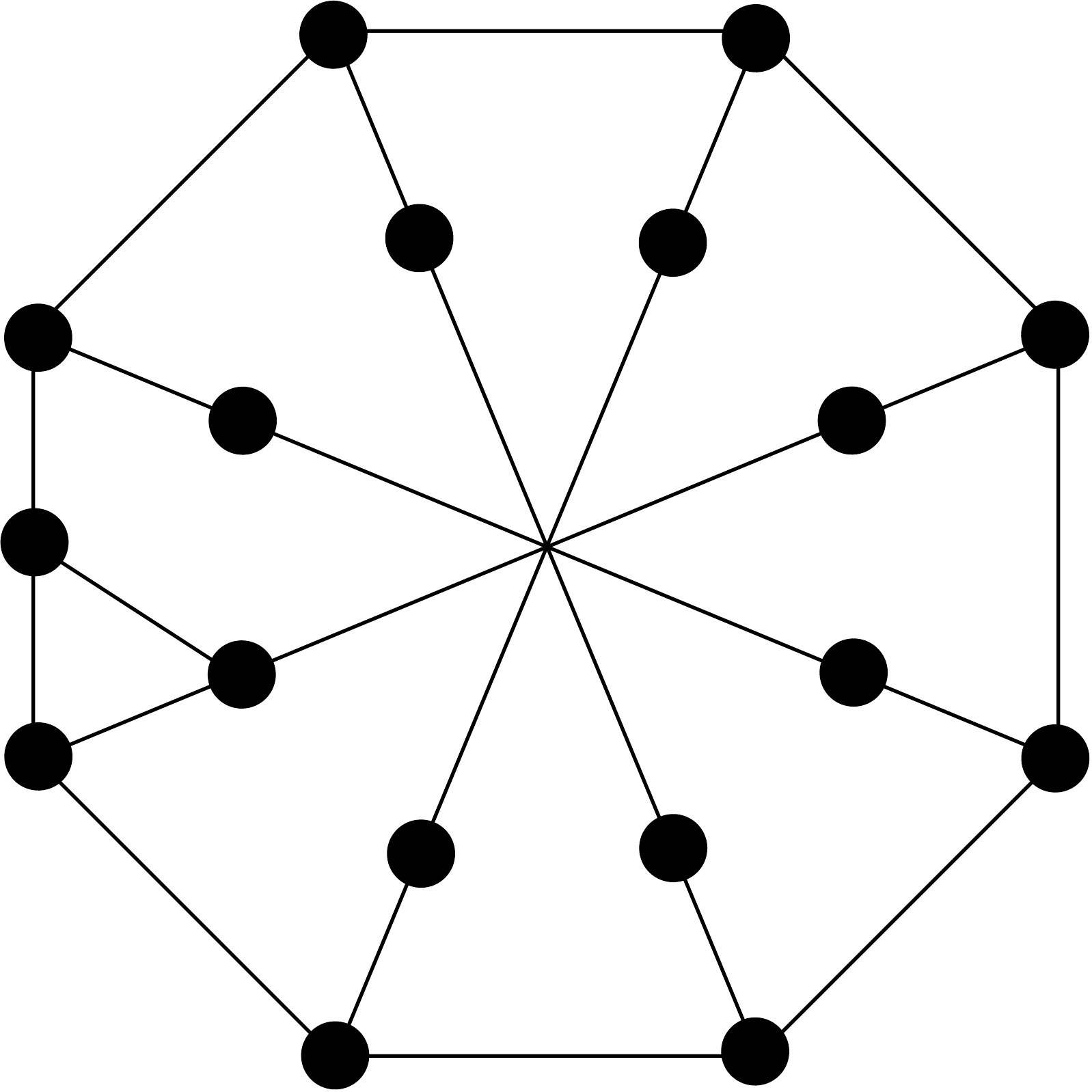} & \includegraphics[height=0.125\textwidth]{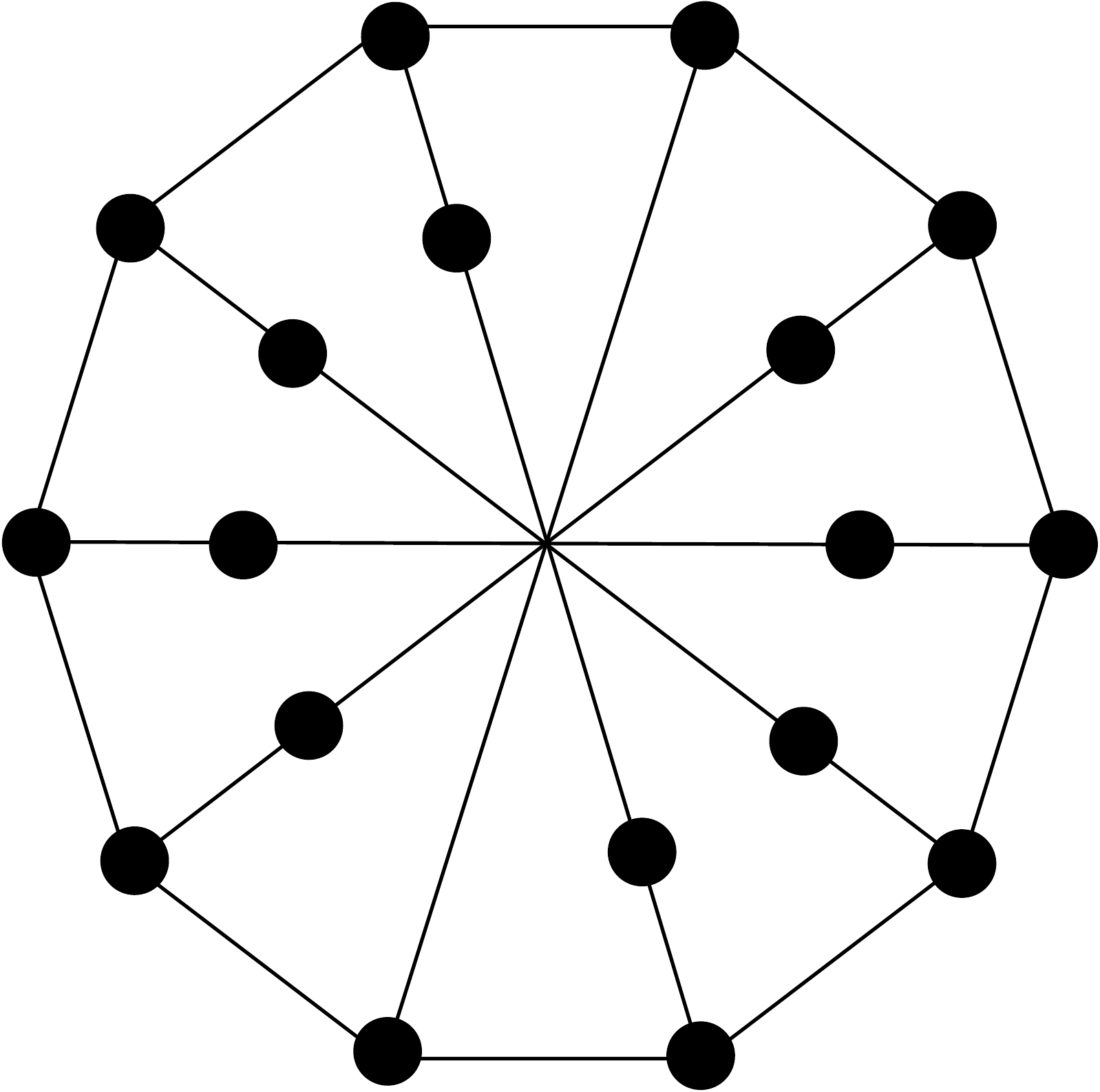} & \includegraphics[height=0.125\textwidth]{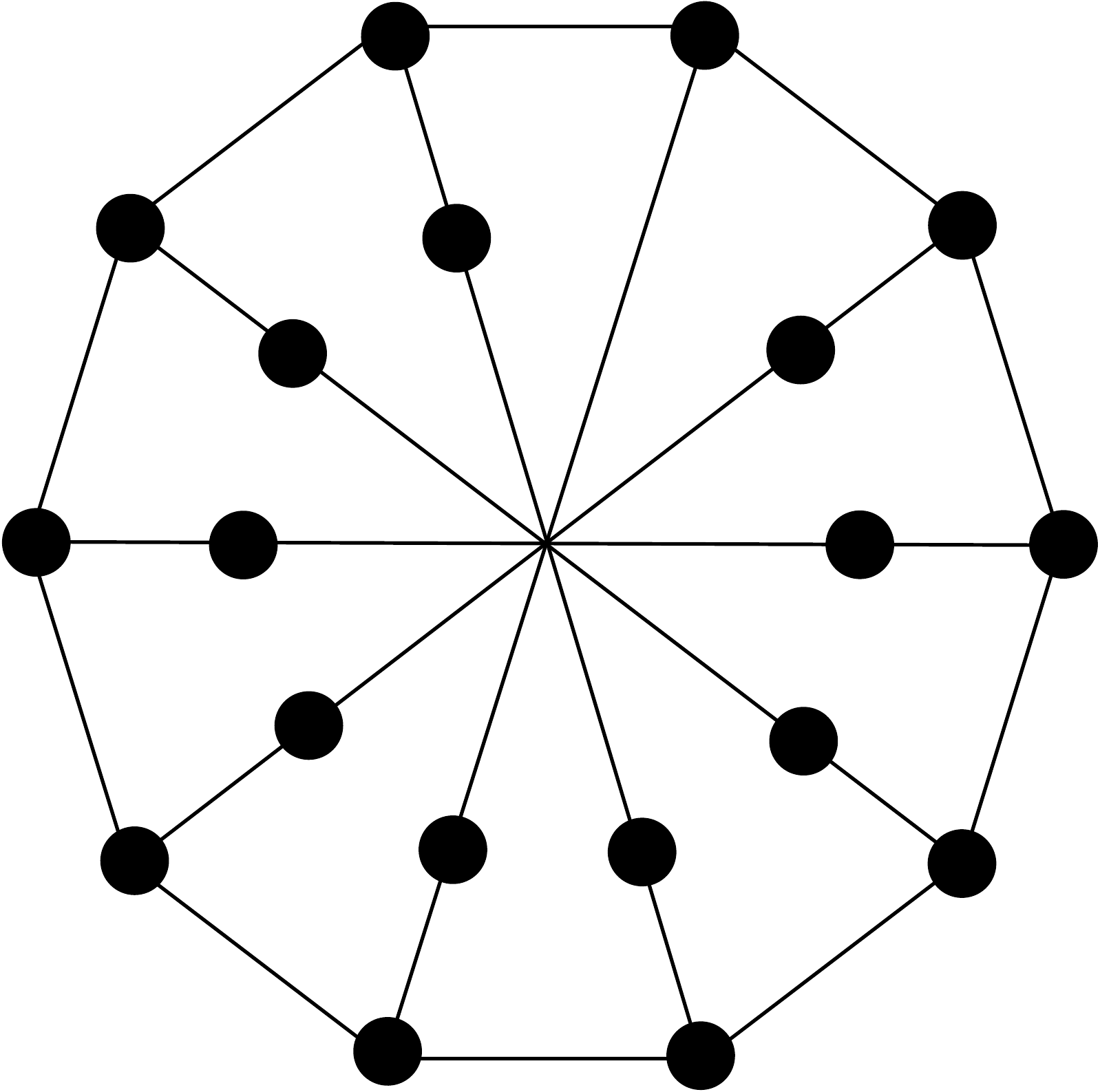} & \includegraphics[height=0.125\textwidth]{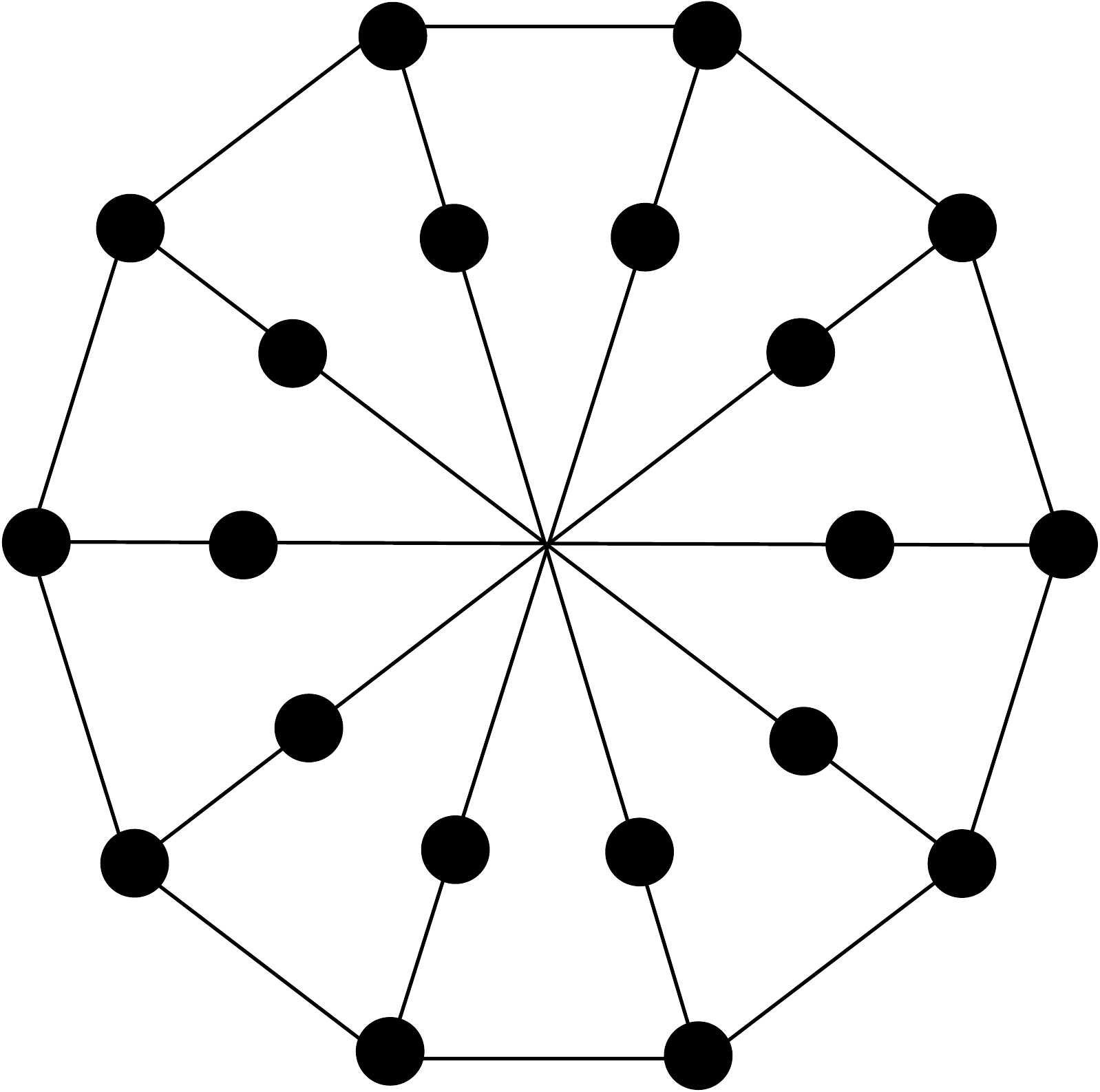} \\
        $n=17$ & $n=18$ & $n=19$ & $n=20$
    \end{tabular}
    \caption{Extremal family of graphs with $\textrm{DET:OLD}(G) = n$ with the smallest number of edges.}
    \label{fig:fam-min-edges}
\end{figure}

The lower bound given in Theorem~\ref{theo:min-edges} on the minimum number of edges in a graph with DET:OLD is sharp for all $n \ge 9$ and Figure~\ref{fig:fam-min-edges} shows a construction for an infinite family of graphs achieving the extremal value for  $9 \le n \le 20$.

\section{NP-completeness of Error-detecting OLD}\label{sec:npc}

It has been shown that many graphical parameters related to detection systems, such as finding the cardinality of the smallest IC, LD, or OLD sets, are NP-complete problems \cite{ld-ic-np-complete-2, NP-complete-ic, NP-complete-ld, old}.
We will now prove that the problem of determining the smallest DET:OLD set is also NP-complete.
For additional information about NP-completeness, see Garey and Johnson \cite{np-complete-bible}.

\npcompleteproblem{3-SAT}{Let $X$ be a set of $N$ variables.
Let $\psi$ be a conjunction of $M$ clauses, where each clause is a disjunction of three literals from distinct variables of $X$.}{Is there is an assignment of values to $X$ such that $\psi$ is true?}

\npcompleteproblem{Error-Detecting Open-locating dominating set (DET-OLD)}{A graph $G$ and integer $K$ with $7 \le K \le |V(G)|$.}{Is there a DET:OLD set $S$ with $|S| \le K$? Or equivalently, is $\textrm{DET:OLD}(G) \le K$?}

\begin{figure}[ht]
    \centering
    \includegraphics[width=0.1\textwidth]{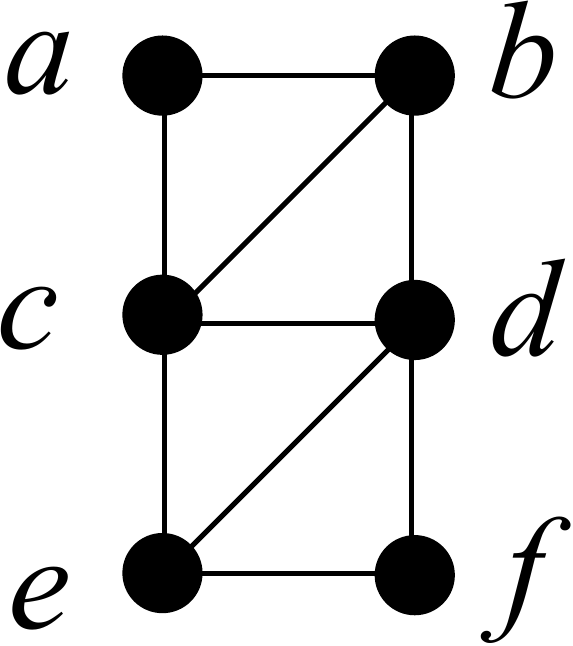}
    \caption{Subgraph $G_6$ }
    \label{fig:g6}
\end{figure}

\begin{lemma}\label{lem:g6}
In the $G_6$ subgraph given in Figure~\ref{fig:g6}, all six vertices internal to $G_6$, as well as at least one external vertex adjacent to $b$ must be detectors in order for DET:OLD to exist.
\end{lemma}
\begin{proof}
Note that vertices $a$, $b$, and $d$ are permitted to have any number of edges going to vertices outside of the $G_6$ subgraph under consideration, but no other edges incident with these 6 vertices are allowed.
Suppose $S$ is a DET:OLD for $G$.
Firstly, to 2-dominate $f$, we require $\{e,d\} \subseteq S$.
To distinguish $e$ and $f$, we require $\{c,f\} \subseteq S$.
To distinguish $c$ and $f$, we require $\{a,b\} \subseteq S$.
Finally, $e$ and $b$ are not distinguished unless $b$ is adjacent to at least one detector external to this $G_6$ subgraph, completing the proof.
\end{proof}

\begin{theorem}
The DET-OLD problem is NP-complete.
\end{theorem}

\begin{wrapfigure}{r}{0.365\textwidth}
    \centering
    \includegraphics[width=0.35\textwidth]{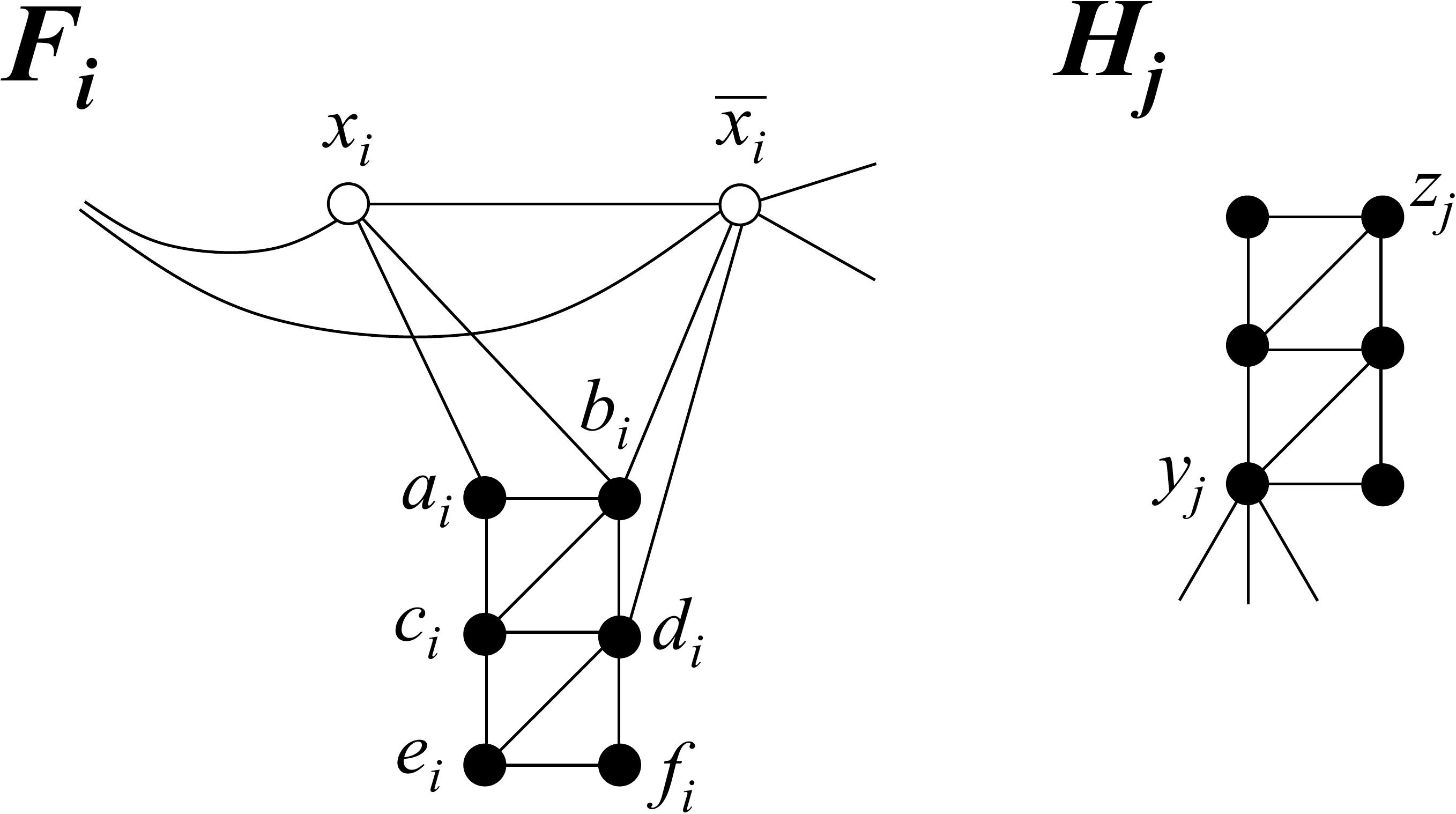}
    \caption{Variable and clause graphs}
    \label{fig:variable-clause}
\end{wrapfigure}
\cbeginproof
Clearly, DET-OLD is NP, as every possible candidate solution can be generated nondeterministically in polynomial time (specifically, $O(n)$ time), and each candidate can be verified in polynomial time using Theorem~\ref{theo:multi-param}.
To complete the proof, we will now show a reduction from 3-SAT to DET-OLD.

\begin{figure}[ht]
    \centering
    \includegraphics[width=0.75\textwidth]{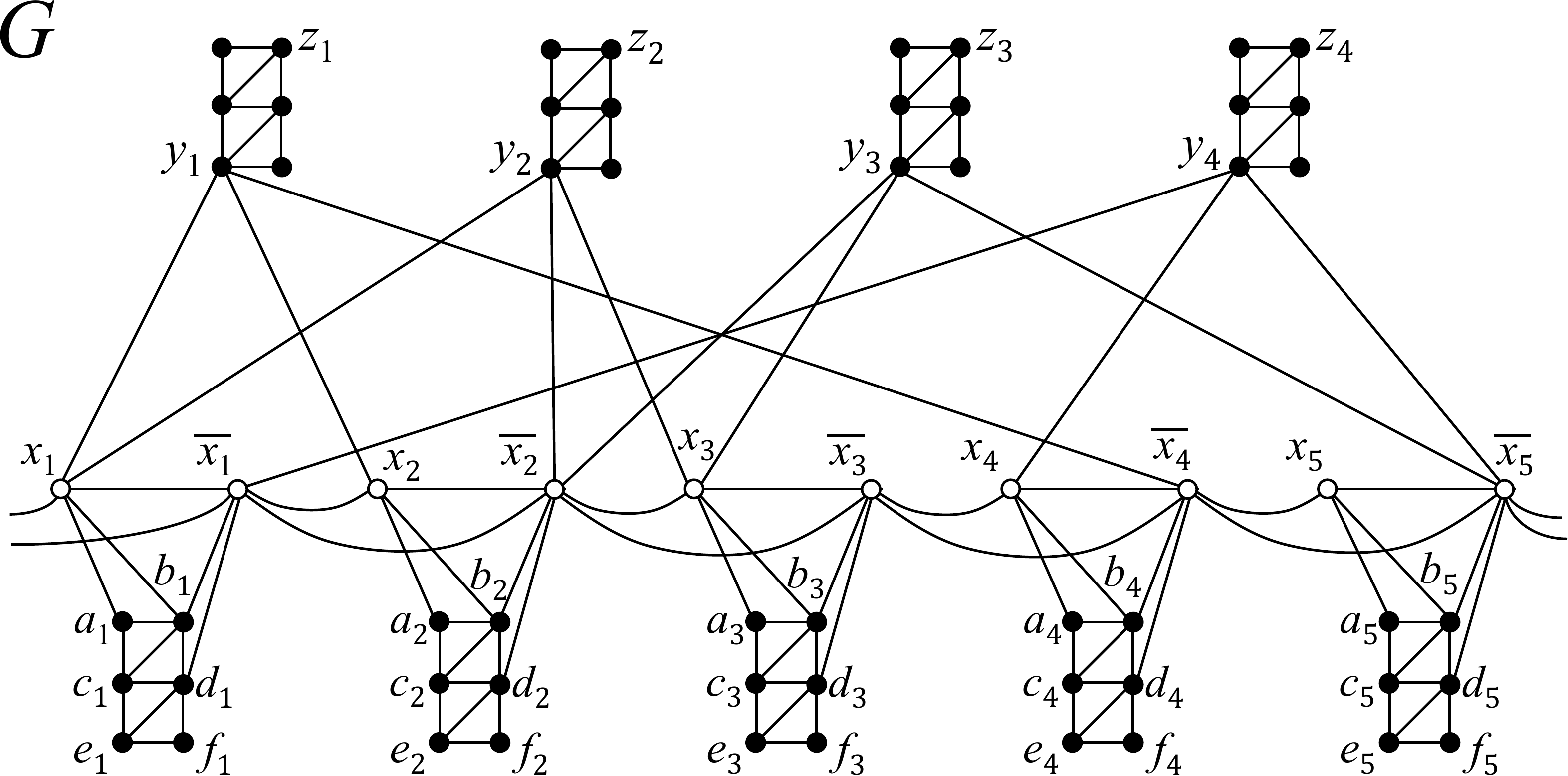}
    \caption{\begin{tabular}[t]{l} Construction of $G$ from $(x_1 \lor x_2 \lor \overline{x_4}) \land (x_1 \lor \overline{x_2} \lor x_3) \land (\overline{x_2} \lor x_3 \lor \overline{x_5}) \land (\overline{x_1} \lor x_4 \lor \overline{x_5})$ \protect\\ with $N = 5$, $M = 4$, $K = 59$ \end{tabular}}
    \label{fig:example-clause}
\end{figure}

Let $\psi$ be an instance of the 3-SAT problem with $M$ clauses on $N$ variables.
We will construct a graph, $G$, as follows.
For each variable $x_i$, create an instance of the $F_i$ graph (Figure~\ref{fig:variable-clause}); this includes a vertex for $x_i$ and its negation $\overline{x_i}$.
For each $i \in \{1, \hdots, N\}$, let $\{\overline{x_i}\;\!x_k, \overline{x_i}\,\overline{x_k}\} \subseteq E(G)$ for $k = (i \!\!\mod N) + 1$.
For each clause $c_j$ of $\psi$, create a new instance of the $H_j$ graph (Figure~\ref{fig:variable-clause}).
For each clause $c_j = \alpha \lor \beta \lor \gamma$, create an edge from the $y_j$ vertex to $\alpha$, $\beta$, and $\gamma$ from the variable graphs, each of which is either some $x_i$ or $\overline{x_i}$; for an example, see Figure~\ref{fig:example-clause}.
The resulting graph has precisely $8N + 6M$ vertices and $16N + 12M$ edges, and can be constructed in polynomial time.
To complete the problem instance, we define $K = 7N + 6M$.

Suppose $S \subseteq V(G)$ is a DET:OLD on $G$ with $|S| \le K$.
By Lemma~\ref{lem:g6}, we require at least the $6N + 6M$ detectors shown by the shaded vertices in Figure~\ref{fig:variable-clause}.
Additionally, Lemma~\ref{lem:g6} gives us the additional requirement that each $b_i$ vertex must be dominated by at least one additional detector outside of its $G_6$ subgraph; thus, for each $i$, $\{x_i,\overline{x_i}\} \cap S \neq \varnothing$, giving us at least $N$ additional detectors.
Thus, we find that $|S| \ge 7N + 6M = K$, implying that $|S| = K$, so $|\{x_i,\overline{x_i}\} \cap S| = 1$ for each $i \in \{1, \hdots, N \}$.
Applying Lemma~\ref{lem:g6} again to the clause graphs yields that each $y_j$ vertex must be dominated by at least one additional detector outside of its $G_6$ subgraph.
As no more detectors may be added, it must be that each $y_j$ is now dominated by one of its three neighbors in the $F_i$ graphs; therefore, $\psi$ is satisfiable.

For the converse, suppose we have a solution to the 3-SAT problem $\psi$; we will show that there is a DET:OLD, $S$, on $G$ with $|S| \le K$.
We construct $S$ by first including all of the $6N + 6M$ vertices as shown in Figure~\ref{fig:variable-clause}.
Next, for each variable, $x_i$, if $x_i$ is true then we let the vertex $x_i \in S$; otherwise, we let $\overline{x_i} \in S$.
Thus, the fully-constructed $S$ has $|S| = 7N + 6M = K$.
Each $b_i$ has its required external dominator due to having $x_i \in S$ or $\overline{x_i} \in S$.
Additionally, because $S$ was constructed from a satisfying truth assignment for the 3-SAT problem, by hypothesis each $y_j$ vertex must also be dominated by one of its (external) term vertices in the $F_i$ subgraphs.
Therefore, each $G_6$ subgraph in $G$ satisfies Lemma~\ref{lem:g6}, and so are internally sufficiently dominated and distinguished.
Since all $G_6$ subgraphs are sufficiently far apart, it is also the case that all vertex pairs in distinct $G_6$ subgraphs are distinguished.
Indeed, it can be shown that all vertices are 2-dominated and $2^\#{}$-distinguished, so $S$ is a DET:OLD set for $G$ with $S \le K$, completing the proof.
\cendproof

\begin{figure}[p]
    \centering
    \begin{tabular}{c@{\hskip 4em}c}
        \centered{\includegraphics[width=0.325\textwidth]{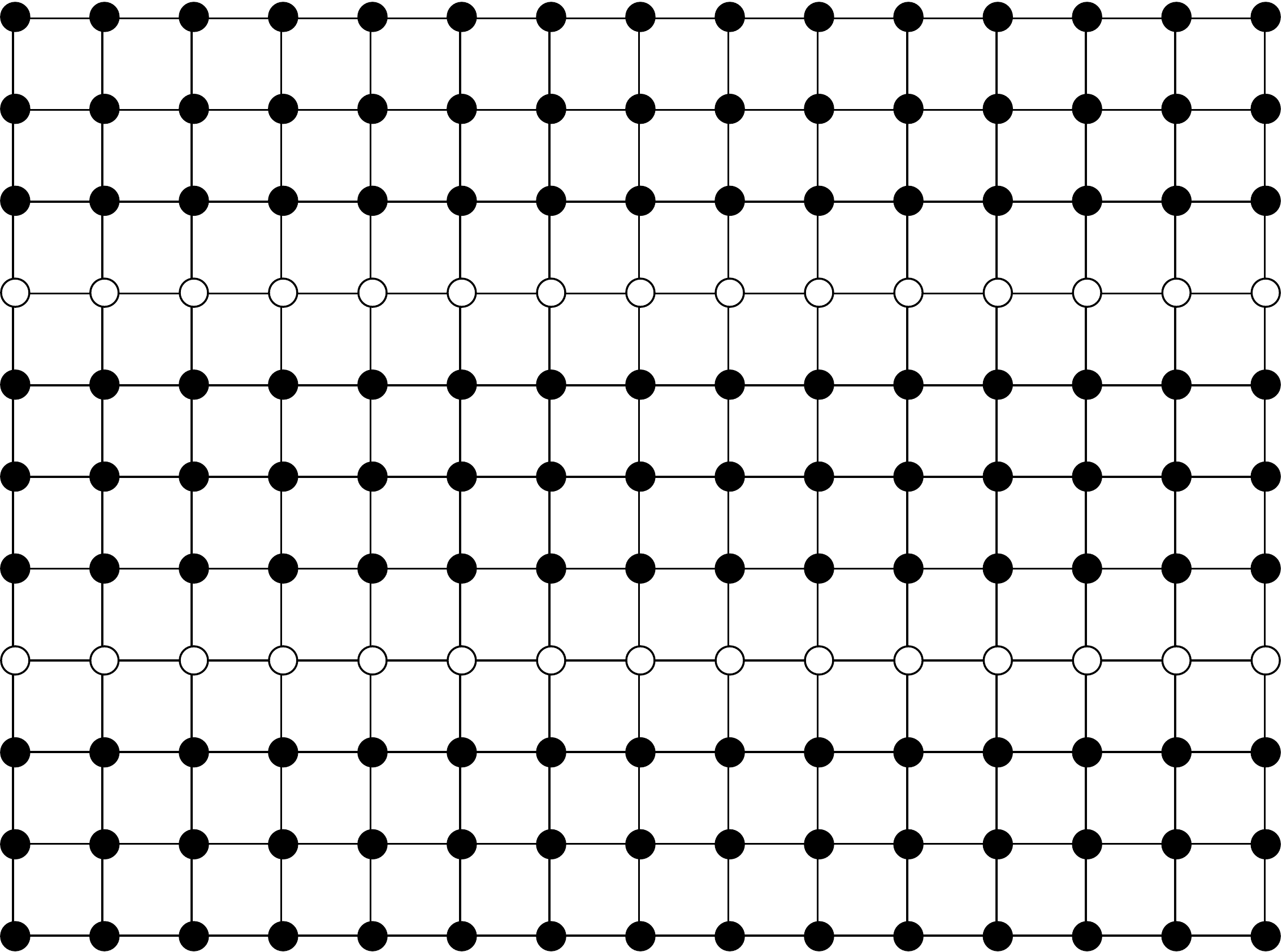}} & \centered{\includegraphics[width=0.325\textwidth]{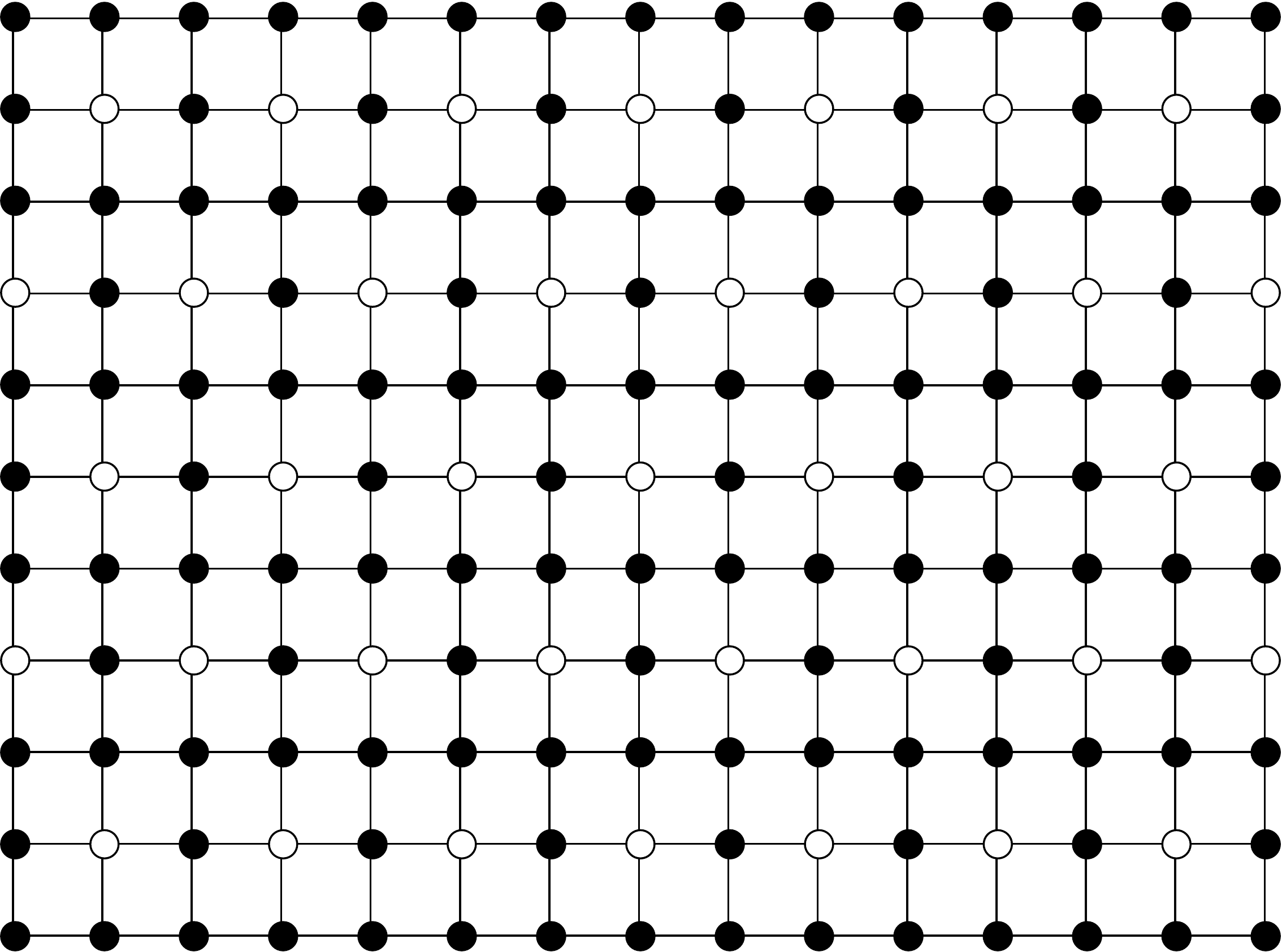}} \\
        (a) & (b) \\ \\
        \centered{\includegraphics[width=0.4\textwidth]{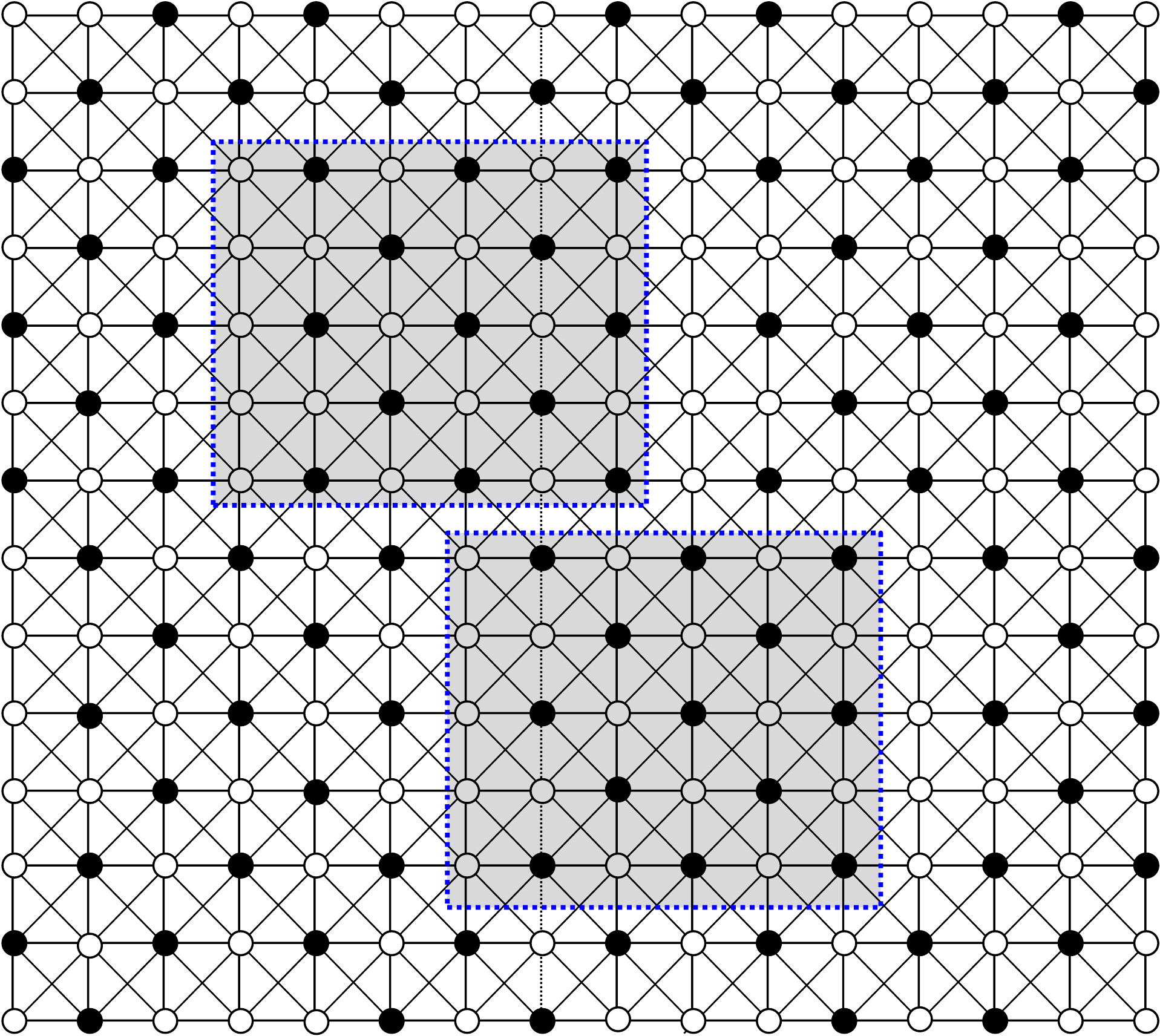}} & \centered{\includegraphics[width=0.4875\textwidth]{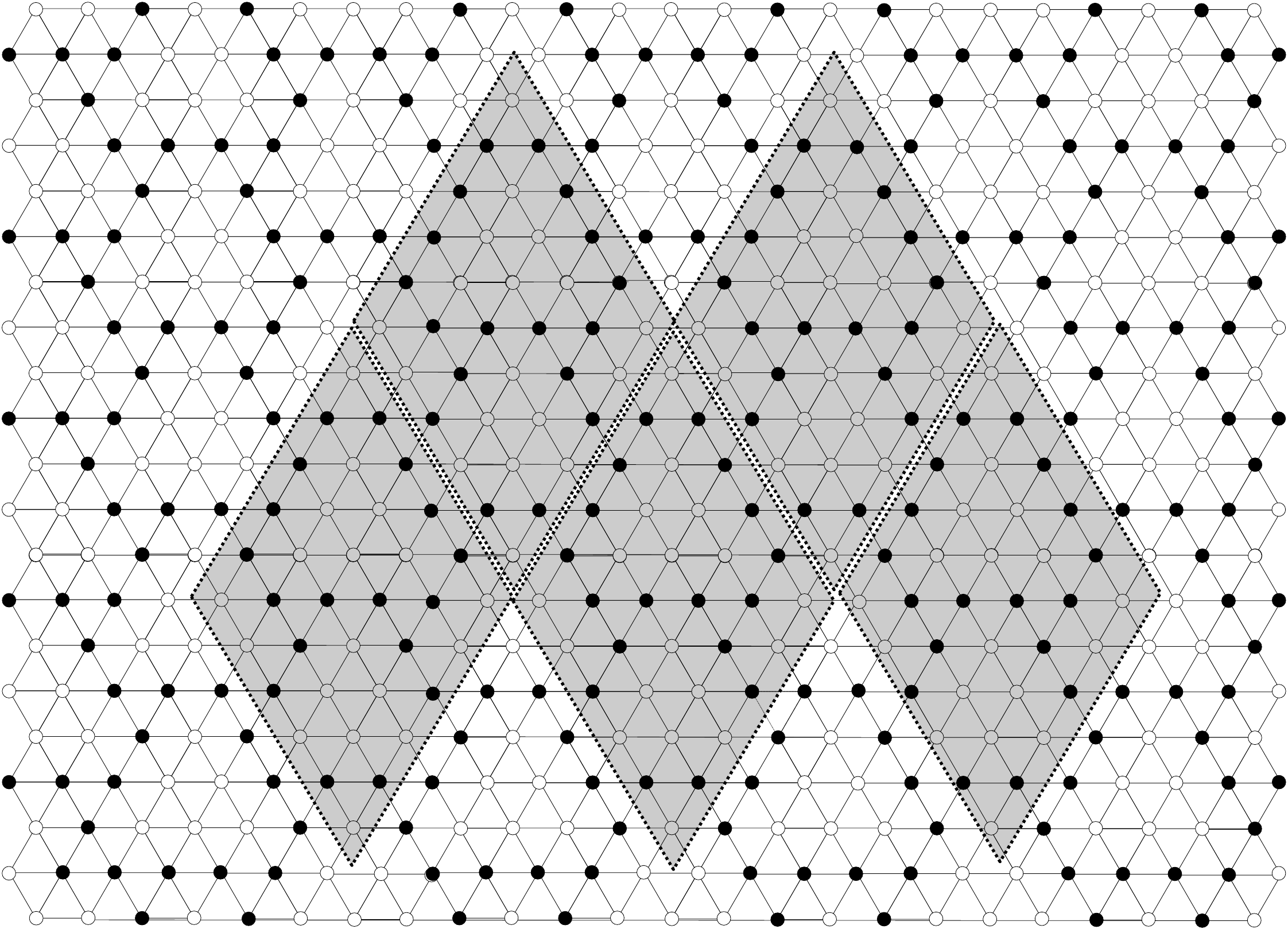}} \\
        (c) & (d) \\ \\
        \centered{\includegraphics[width=0.375\textwidth]{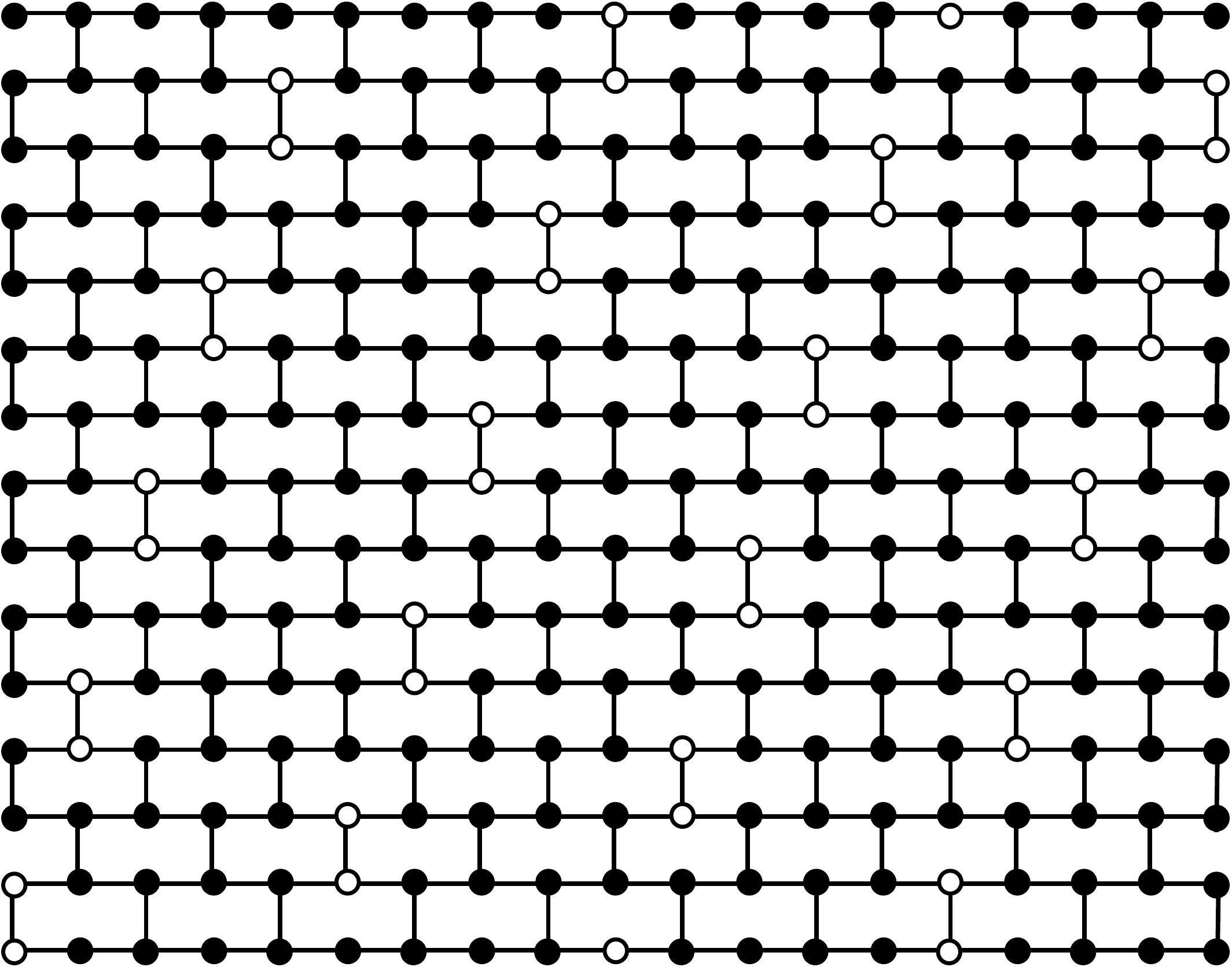}} & \centered{\includegraphics[width=0.375\textwidth]{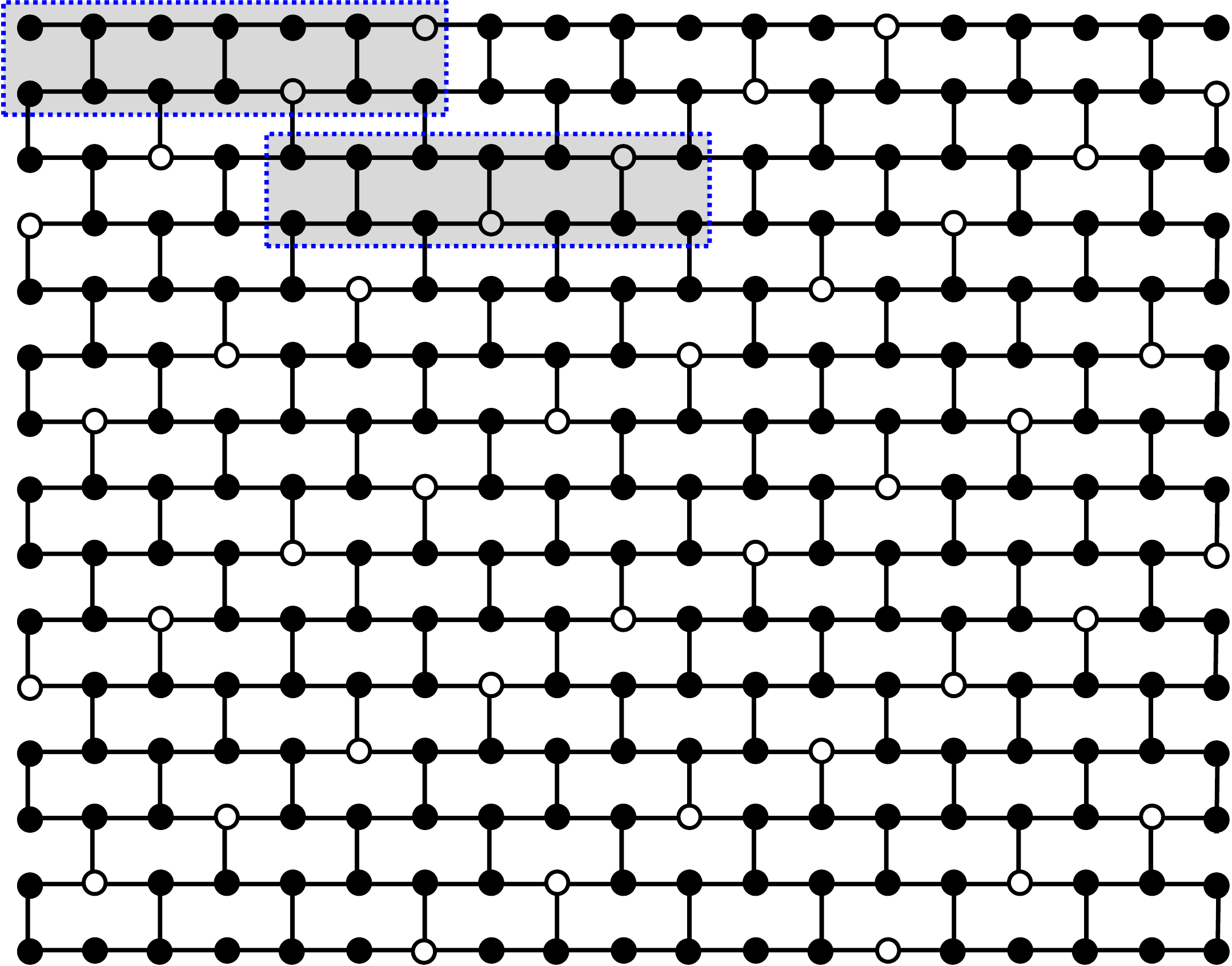}} \\
        (e) & (f)
    \end{tabular}
    \caption{Our best constructions of DET:OLD on SQR (a) and (b), KNG (c), TRI (d), and HEX (e) and (f). Shaded vertices denote detectors.}
    \label{fig:inf-grids-det-old-solns}
\end{figure}

\section{DET:OLD sets for Infinite Grids}\label{sec:grids}

Theorem~\ref{theo:det-old-uppers} presents upper bounds and some tight bounds for DET:OLD on several infinite grids; solutions achieving these upper bounds are presented in Figure~\ref{fig:inf-grids-det-old-solns}.
Solutions (b), (c), and (f) are new solutions, while the others were first presented in other papers \cite{ourtri, ftold}.
The solutions for the infinite hexagonal grid (HEX), the infinite square grid (SQR), and the infinite triangular grid (TRI) are tight bounds, while the exact value for the infinite king grid (KNG) is currently unknown.
Subfigure (c) gives the best (lowest-density) solution we have found for KNG.
We believe that standard discharging share arguments can be used to prove a lower bound of $\frac{60}{151}$ for KNG.

\begin{theorem}\label{theo:det-old-uppers}
The upper and lower bounds on DET:OLD:
\begin{enumerate}[label=\roman*]
    \item For the infinite hexagonal grid HEX, $\textrm{DET:OLD\%}(HEX) = \frac{6}{7}$. (Seo and Slater \cite{ftold})
    \item For the infinite square grid SQR, $\textrm{DET:OLD\%}(SQR) = \frac{3}{4}$. (Seo and Slater \cite{ftold})
    \item For the infinite triangular grid TRI, $\textrm{DET:OLD\%}(TRI) = \frac{1}{2}$. (Jean and Seo \cite{ourtri})
    \item For the infinite king grid KNG,
$\textrm{DET:OLD\%}(KNG) \le \frac{13}{30}$.
\end{enumerate}
\end{theorem}

\section{Extremal cubic graphs}\label{sec:cubic}

In this section, we characterize cubic graphs that permit a DET:OLD set.
We also consider extremal cubic graphs on $\textrm{DET:OLD}(G)$.

\begin{observation}\label{obs:c4-free}
A cubic graph $G$ has DET:OLD if and only if $G$ is $C_4$-free.
\end{observation}
\begin{proof}
Let $abcd$ be a 4-cycle in $G$.
We see that $a$ and $c$ cannot be distinguished, a contradiction.
For the converse, we will show that $S=V(G)$ is a DET:OLD for a $C_4$-free cubic graph.
Let $u,v \in V(G)$; we know that $|N(u) \cap N(v)| \le 1$ because $G$ is $C_4$-free.
Thus, $|N(u) - N(v)| \ge 2$, implying $u$ is distinguished from $v$, completing the proof.
\end{proof}

\begin{figure}[ht]
    \centering
    \begin{tabular}{c@{\hskip 4em}c@{\hskip 4em}c}
        \centered{\includegraphics[width=0.14\textwidth]{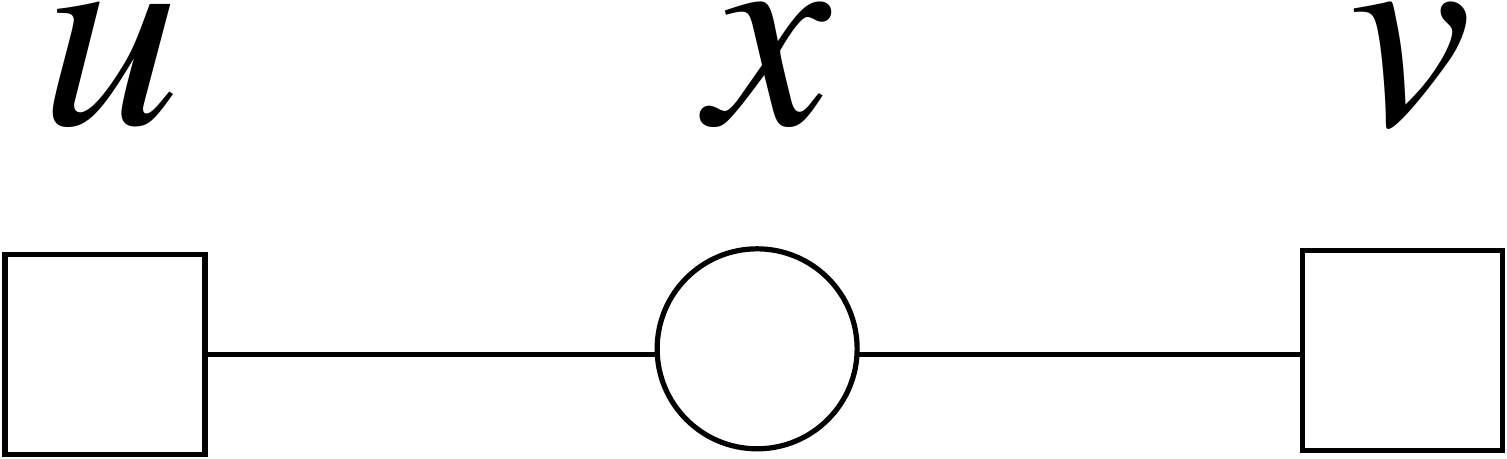}} & \centered{\includegraphics[width=0.25\textwidth]{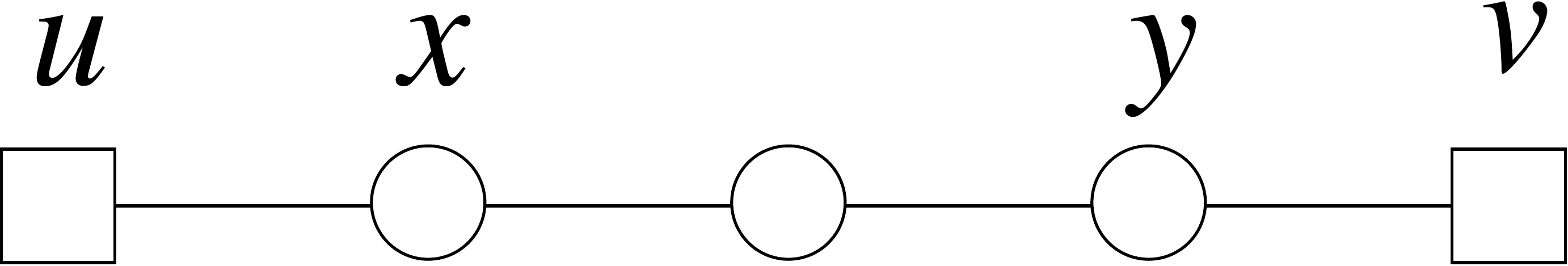}} & \centered{\includegraphics[width=0.13\textwidth]{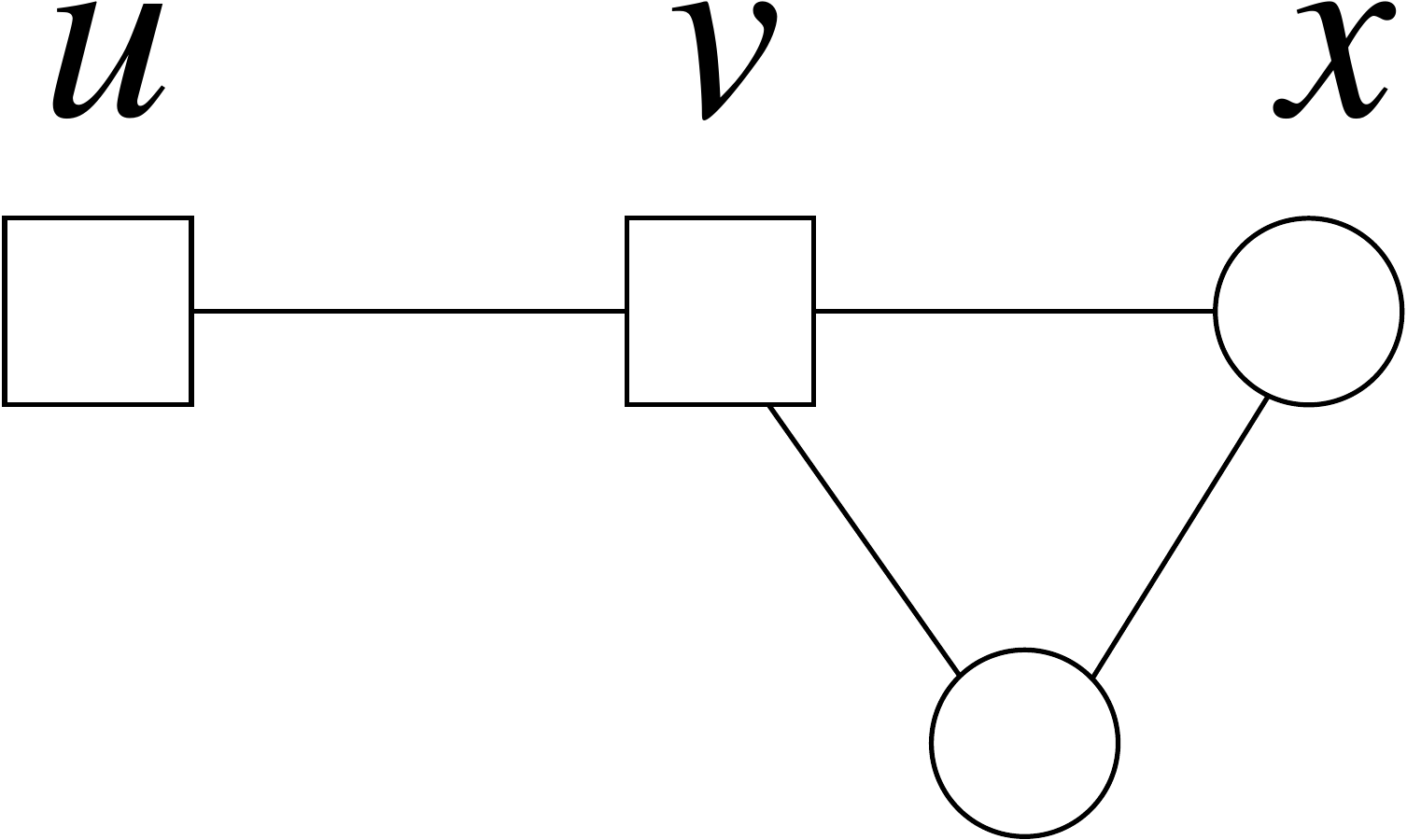}} \\
        (a) & (b) & (c)
    \end{tabular}
    \caption{The three nonisomorphic constructions of length 2 and 4 trails from $u$ to $v$}
    \label{fig:det-old-t24-cases}
\end{figure}

\begin{observation}\label{obs:t24}
Let $S$ be a DET:OLD set on cubic graph $G$.
For all vertices $v \in V(G)-S$, we have $T_2(v) \cup T_4(v) \subseteq S$.
\end{observation}
\begin{proof}
Assume to the contrary that there exist $u,v \in V(G)-S$ such that $u \in T_2(v) \cup T_4(v)$.
We know that $u \neq v$ because $G$ is $C_4$-free due to the existence of DET:OLD.
We consider the three nonisomorphic ways to form a length 2 or 4 trail from $u$ to $v$ as illustrated in Figure~\ref{fig:det-old-t24-cases}.
We see that in Figure~\ref{fig:det-old-t24-cases}~(a), $x$ cannot be 2-dominated, in (b), $x$ and $y$ cannot be distinguished, and in (c), $v$ and $x$ cannot be distinguished.
All three cases contradict the existence of $S$, completing the proof.
\end{proof}

\begin{theorem}\label{theo:t24}
Let $G$ be a $C_4$-free cubic graph, and let $\overline{S} \subseteq V(G)$ such that for all distinct $u,v \in \overline{S}$, $u \notin T_2(v) \cup T_4(v)$. Then $S = V(G) - \overline{S}$ is a DET:OLD on $G$.
\end{theorem}
\begin{proof}
Let $u \in V(G)$.
We see that for any distinct $x,y \in N(u)$, $x \in T_2(y)$.
Thus, by assumption it must be that $|N(u) \cap \overline{S}| \le 1$, implying that $dom(u) \ge 2$.
We now know that all vertices are at least 2-dominated by $S$.
Next, we consider the following three cases depending on the distance between a pair of vertices and show the pair is $2^{\#}$-distinguished.

\textbf{Case 1:} Suppose $u,v \in V(G)$ with $d(u,v) = 1$.
Because $G$ is twin-free (due to being a $C_4$-free cubic graph) and regular, it must be that $\exists x \in N(u) - N[v]$ and $\exists y \in N(v) - N[u]$.
Suppose $\exists w \in N(u) \cap N(v)$.
Then we see that $\forall p \in N[u] \cup N[v]$, $(N[u] \cup N[v]) - \{p\} \subseteq T_2(p) \cup T_4(p)$, implying that $|(N[u] \cup N[v]) \cap \overline{S}| \le 1$, from which it can be seen that $u$ and $v$ must be distinguished.
Otherwise, we can assume that $N(u) \cap N(v) = \varnothing$.
If $u \in \overline{S}$, then $N(v) - N[u] \subseteq T_2(u) \subseteq S$, so $v$ is distinguished for $u$ and we would be done; otherwise we assume $u \in S$ and by symmetry $v \in S$.
We see that $|(N(u) - N[v]) \cap \overline{S}| \le 1$, so $u$ is distinguished from $v$.

\textbf{Case 2:} Suppose $u,v \in V(G)$ with $d(u,v) = 2$, and let $N(u) = N(v) = \{w\}$, as $G$ is $C_4$-free.
Let $N(u) - N[v] = \{u',u''\}$ and let $N(v) - N[u] = \{v', v''\}$.
We see that $\forall p \in \{u',u'',v',v''\}$, $\{u',u'',v',v''\} - \{p\} \subseteq T_2(p) \cup T_4(p) \subseteq S$, which implies that $|\{u',u'',v',v''\} \cap \overline{S}| \le 1$.
Therefore, $u$ and $v$ must be distinguished.

\textbf{Case 3:} Suppose $u,v \in V(G)$ with $d(u,v) \ge 3$, and let $N(u) = \{u',u'',u'''\}$ and $N(v) = \{v',v'',v'''\}$.
We see that for any $p \in N(u)$, $\{u',u'',u'''\} - \{p\} \subseteq T_2(p) \subseteq S$.
Thus $u$ and $v$ must be distinguished, completing the proof.
\end{proof}

From Observations \ref{obs:c4-free} and \ref{obs:t24} and Theorem~\ref{theo:t24}, we have the following full characterization of DET:OLD in cubic graphs.

\begin{corollary}\label{cor:tchar}
Let $G$ be a cubic graph and $S \subseteq V(G)$.
Then $S$ is a DET:OLD if and only if $G$ is $C_4$-free and for all distinct $u,v \in V(G)-S$, $u \notin T_2(v) \cup T_4(v)$.
\end{corollary}

\medskip
The upper bound on $\textrm{DET:OLD}(G)$ for cubic graphs is known to be  $\frac{45}{46}n$ \cite{ft-old-cubic}, and we can improve it using Corollary~\ref{cor:tchar}.

\begin{theorem}
If $G$ is a cubic graph that permits DET:OLD, then $\textrm{DET:OLD\%}(G) \le \frac{30}{31}$.
\end{theorem}
\begin{proof}
We will show that we can construct a set $\overline{S}$ with the property that $\forall v \in V(G)$, $\exists u \in \overline{S}$ such that $v \in T_0(u) \cup T_2(u) \cup T_4(u)$.
Because $|T_0(u)| = 1$,  $|T_2(u)| \le 6$, and $|T_4(u)| \le 24$, this construction will result in a detector set $S = V(G) - \overline{S}$ with density at most $\frac{6 + 24}{1 + 6 + 24} = \frac{30}{31}$.
Assume to the contrary that we have a maximal $\overline{S}$ set, but $\exists x \in V(G)$ such that $x \notin T_0(u) \cup T_2(u) \cup T_4(u)$ for any $u \in \overline{S}$, implying that $x \notin \overline{S}$.
Then by Corollary~\ref{cor:tchar}, we see that $\overline{S} \cup \{x\}$ still satisfies the requirements of our $\overline{S}$ set, contradicting maximality of $\overline{S}$.
Therefore, we have a DET:OLD set $S$ with density at most $\frac{30}{31}$, completing the proof.
\end{proof}

\begin{corollary}\label{cor:det-old-n-1}
If $G$ is a cubic graph that permits DET:OLD, then $\textrm{DET:OLD}(G) \le n - 1$.
\end{corollary}

\begin{figure}[ht]
    \centering
    \begin{tabular}{c@{\hskip 0.9em}c@{\hskip 0.9em}c@{\hskip 0.9em}c@{\hskip 0.9em}c}
        \centered{\includegraphics[width=0.1425\textwidth]{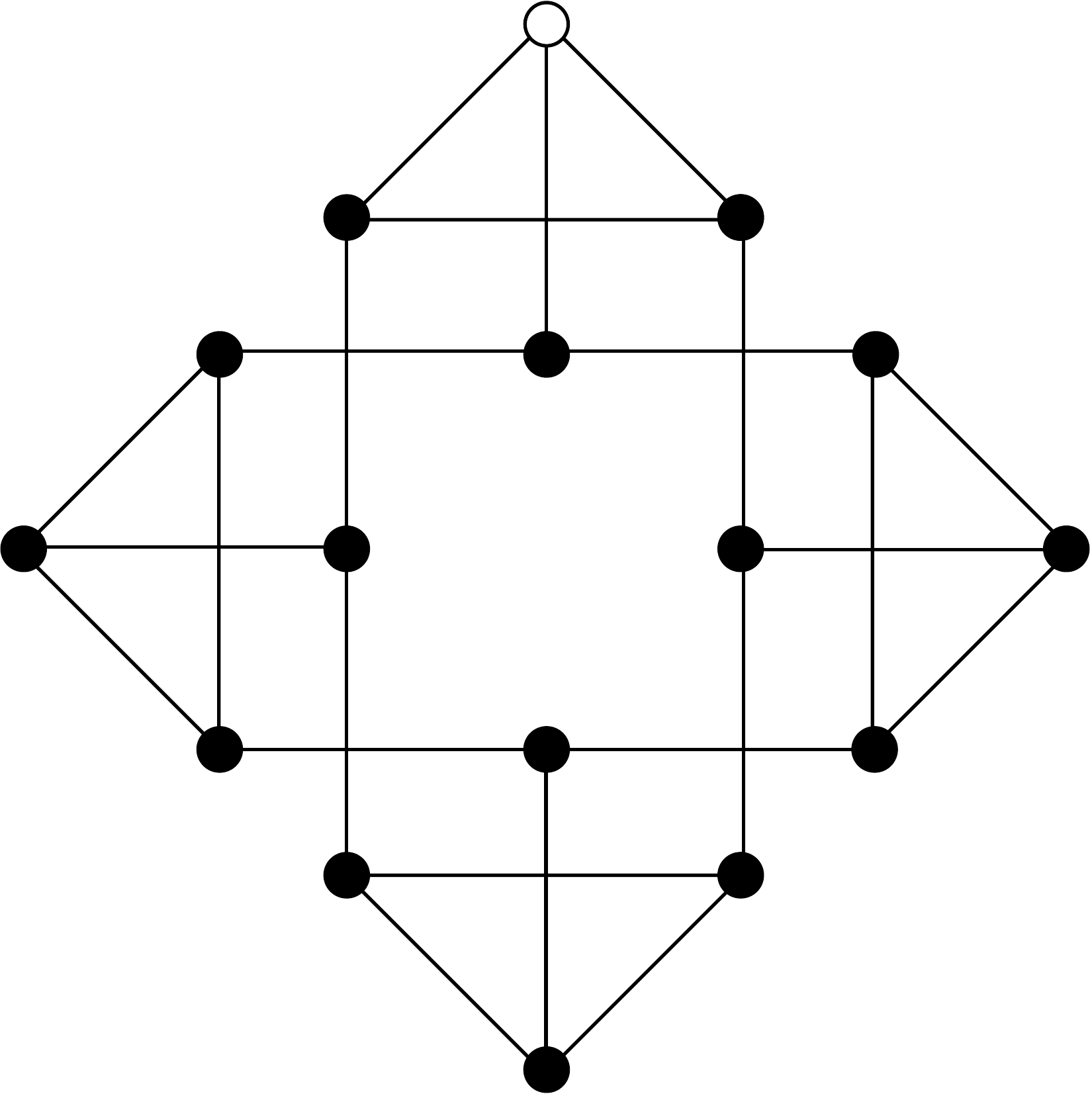}} & \centered{\includegraphics[width=0.1425\textwidth]{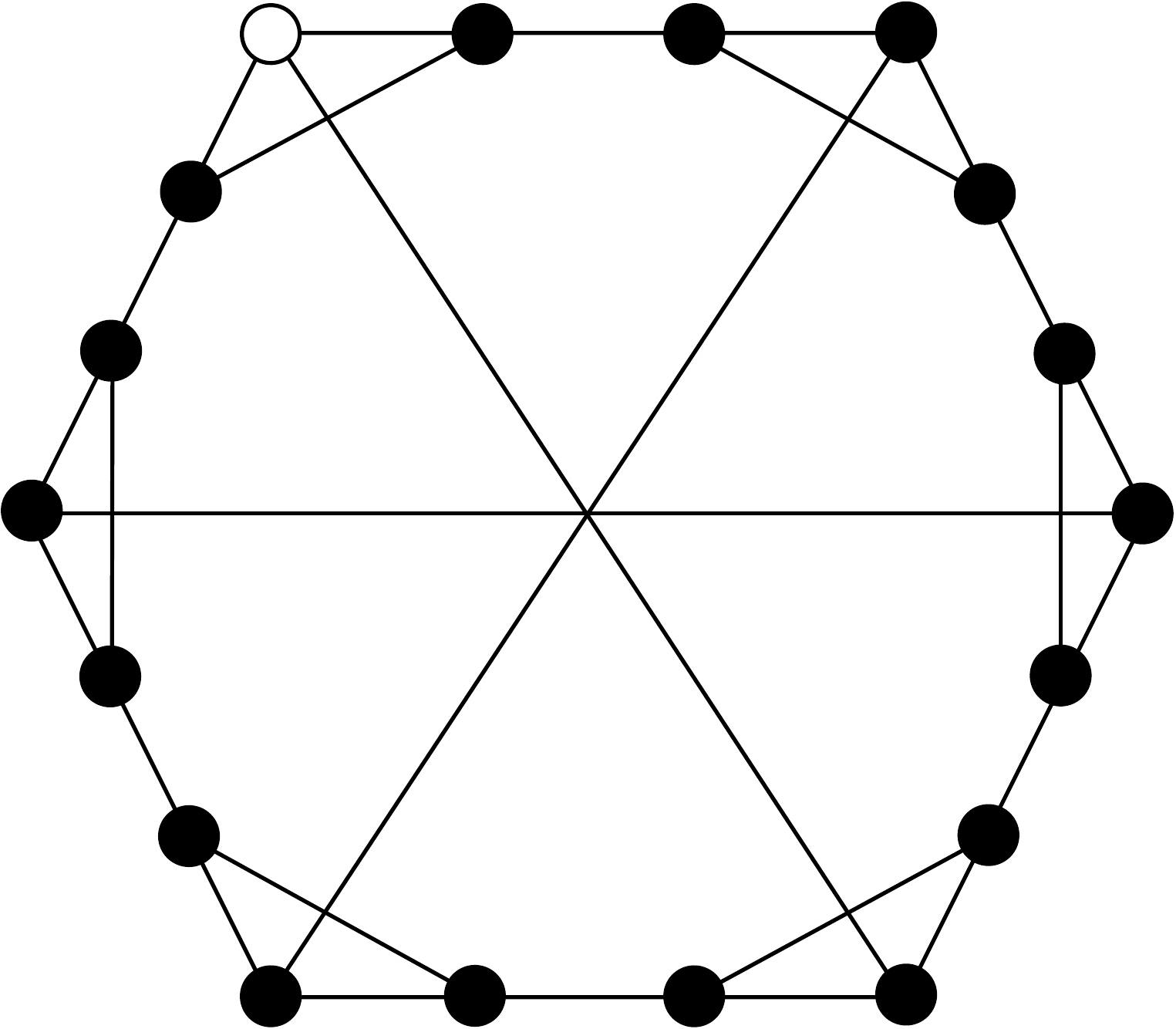}} & \centered{\includegraphics[width=0.165\textwidth]{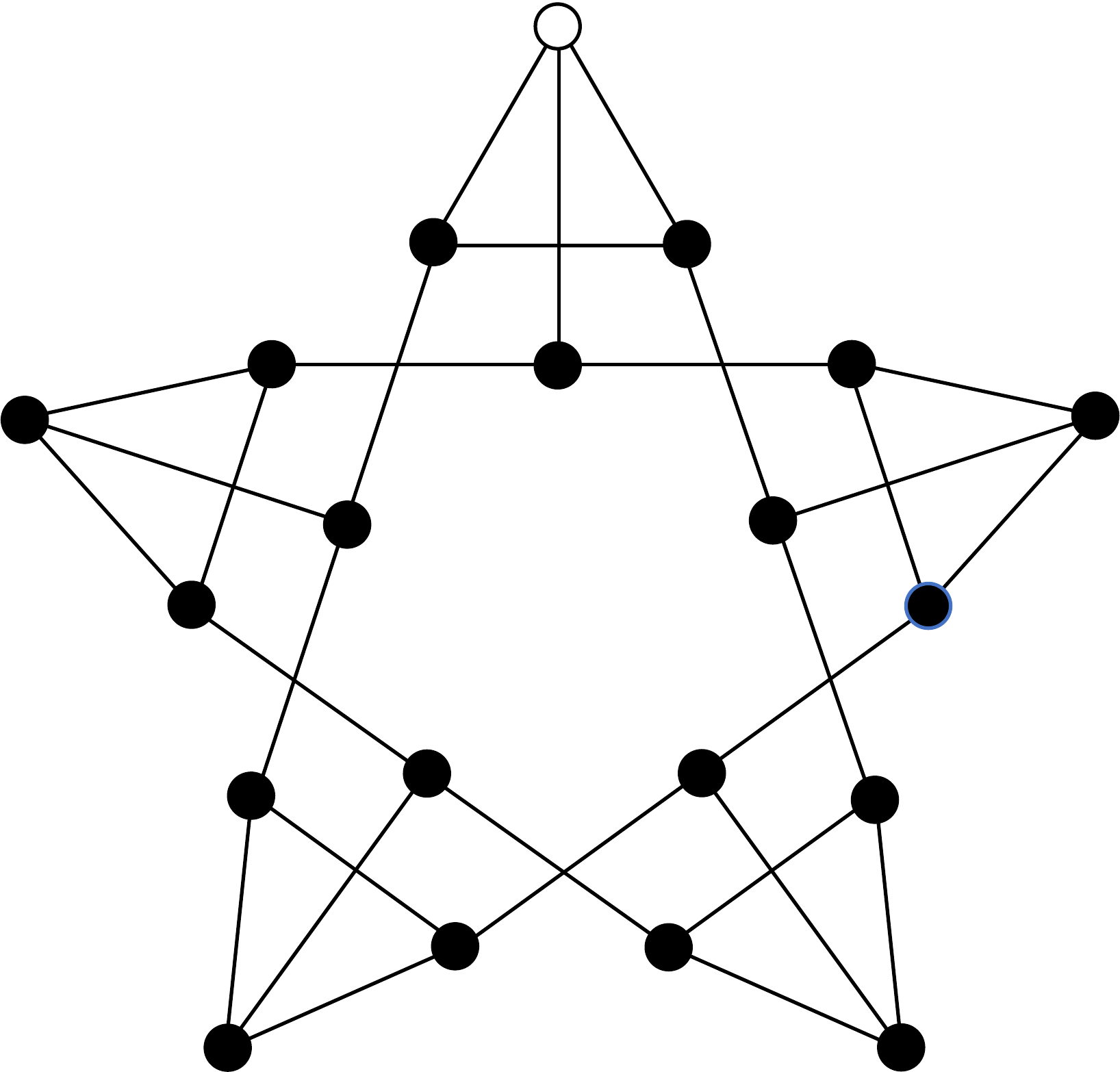}} & \centered{\includegraphics[width=0.16\textwidth]{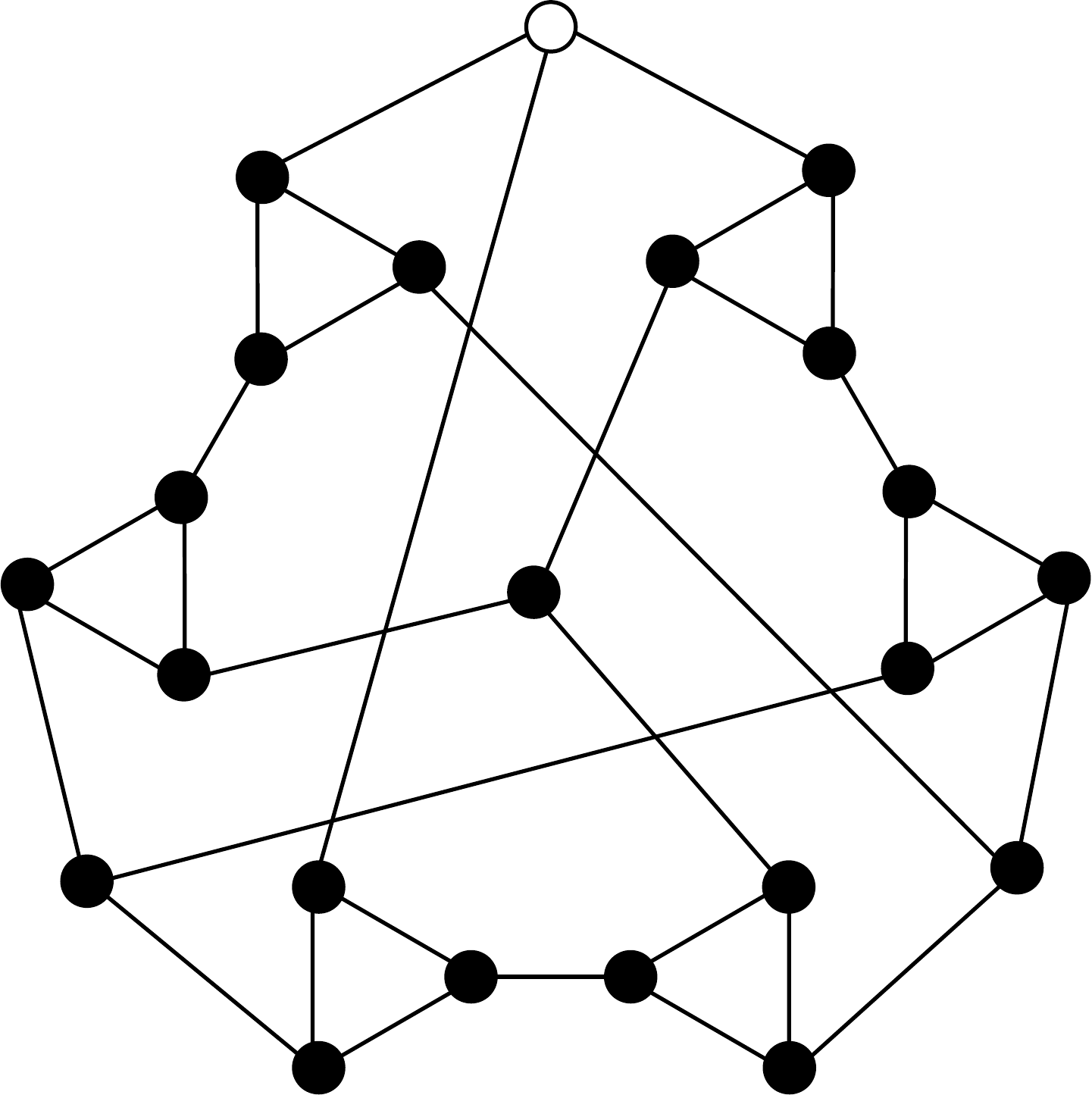}} & \centered{\includegraphics[width=0.16\textwidth]{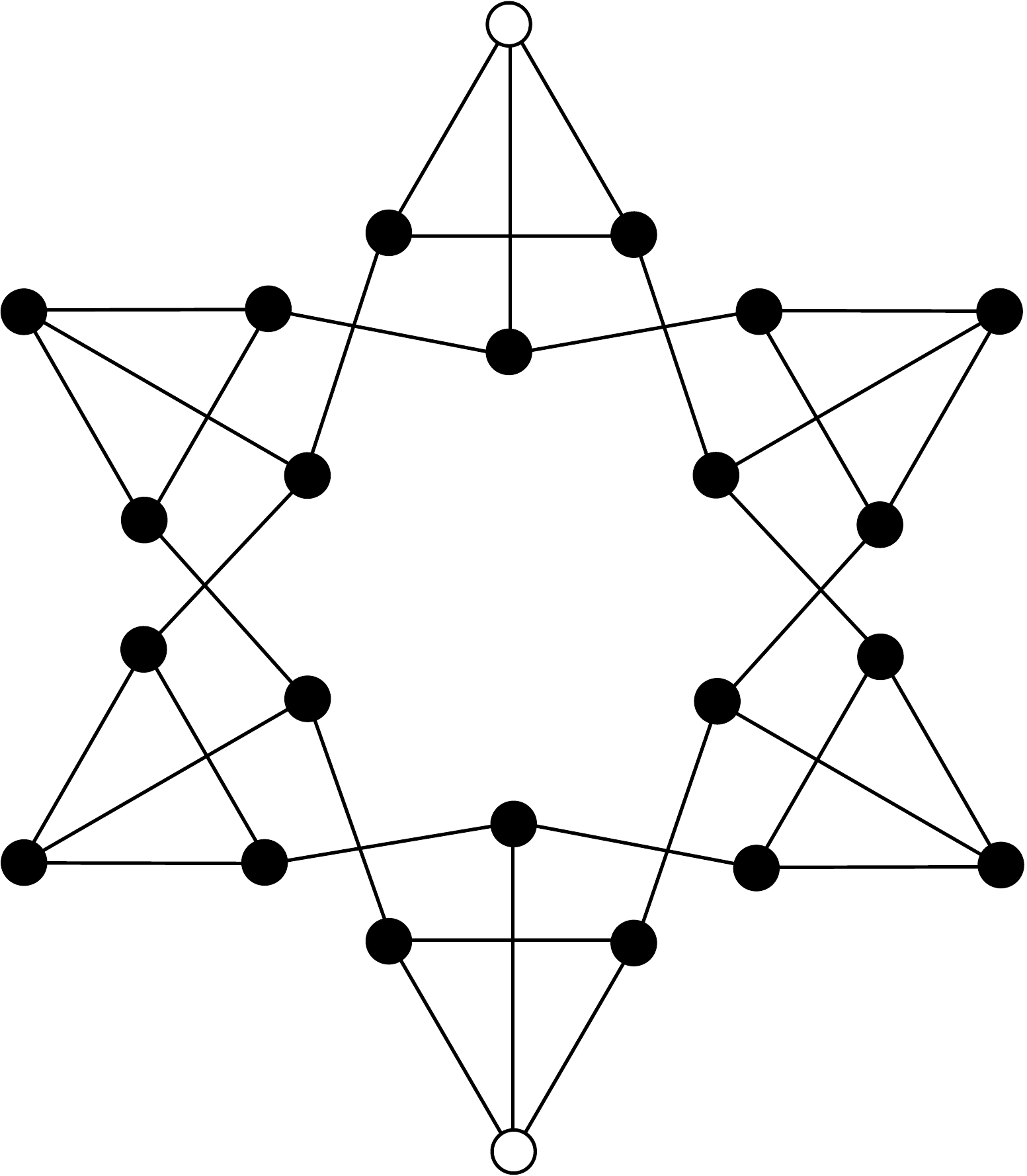}} \\
        $15/16$ & $17/18$ & $19/20$ & $21/22$ & $22/24$
    \end{tabular}
    \caption{Extremal cubic graphs on $n$ vertices with the highest density for $16 \le n \le 24$}
    \label{fig:det-old-ext-fam}
\end{figure}

Figure~\ref{fig:det-old-ext-fam} shows extremal cubic graphs with the highest density on $n$ vertices for $16 \le n \le 24$.
The $n = 22$ graph shown in Figure~\ref{fig:det-old-ext-fam} has the highest density we have found so far of $\frac{21}{22}$, and we conjecture the density $\frac{21}{22}$ is the tight upper bound for all cubic graphs.

\begin{conjecture}
If $G$ is a $C_4$-free cubic graph, then $\textrm{DET:OLD}\%(G) \le \frac{21}{22}$.
\end{conjecture}

\bibliographystyle{ACM-Reference-Format}
\bibliography{refs}

\end{document}